\documentclass[preprint,3p]{elsarticle}
\usepackage[pdftex]{hyperref}
\hypersetup{
	colorlinks   = true, 
	urlcolor     = blue, 
	linkcolor    = blue, 
	citecolor   = red 
}
\usepackage[labelformat=simple]{subcaption}
\usepackage{graphics}
\usepackage{epsfig}

\usepackage{float}
\usepackage{amssymb}
\usepackage{amsmath}
\usepackage{amsthm}
\usepackage{amsfonts}
\usepackage{dsfont}
\usepackage{algorithm}
\usepackage{algorithmic}
\usepackage{xcolor}
\usepackage[noabbrev]{cleveref}
\newtheorem{ssmptn}{Assumption}

\newtheorem{prpstn}{Proposition}
\newtheorem{dfntn}{Definition}
\newtheorem{thrm}{Theorem}
\newtheorem{rmrk}{Remark}

\usepackage{mathrsfs}
\usepackage{enumitem}

\biboptions{sort&compress}

\newcommand{\Vn}{\mathbf{n}}
\newcommand{\Vb}{\mathbf{b}}

\newcommand{\R}{\mathbb{R}}

\newcommand{\Par}{\mathscr{P}}

\newcommand{\Sk}{\mathscr{S}}
\newcommand{\Ak}{\mathscr{A}}
\newcommand{\Vk}{\mathscr{V}}

\renewcommand{\div}{\mathrm{div}}

\newcommand{\llbrace}{\lbrace\hspace{-.12cm}\lbrace }
\newcommand{\rrbrace}{\rbrace \hspace{-.12cm}\rbrace }
\newcommand{\llbracket}{[\![}
\newcommand{\rrbracket}{]\!]}
\DeclareMathOperator*{\argmin}{arg\,min}

\usepackage{numcompress}

\journal{Finite Elements in Analysis and Design}

\begin{document}

\begin{frontmatter}
\title{An automatic-adaptivity stabilized finite element method via
  residual minimization for heterogeneous, anisotropic advection-diffusion-reaction problems}
\author[one,three]{Roberto J. Cier\corref{mycorrespondingauthor}}
\cortext[mycorrespondingauthor]{Corresponding author}
\ead{rcier93@gmail.com}

\author[two]{Sergio Rojas}

\author[two,three]{Victor M. Calo}

\address[one]{School of Civil and Mechanical Engineering, Curtin University, Kent Street, Bentley, Perth, WA 6102, Australia}

\address[two]{School of Earth and Planetary Sciences, Curtin University, Kent Street, Bentley, Perth, WA 6102, Australia}

\address[three]{Mineral Resources, Commonwealth Scientific and Industrial Research Organisation (CSIRO), Kensington, Perth, WA 6152, Australia}

\begin{abstract}

In this paper, we describe a stable finite element formulation for advection-diffusion-reaction problems that allows for robust automatic adaptive strategies to be easily implemented. We consider locally vanishing, heterogeneous, and anisotropic diffusivities, as well as advection-dominated diffusion problems. The general stabilized finite element framework was described and analyzed in~\cite{ calo2019} for linear problems in general, and tested for pure advection problems. The method seeks for the discrete solution through a residual minimization process on a proper stable discontinuous Galerkin (dG) dual norm. This technique leads to a saddle-point problem that delivers a stable discrete solution and a robust error estimate that can drive mesh adaptivity. In this work, we demonstrate the efficiency of the method in extreme scenarios, delivering stable solutions. The quality and performance of the solutions are comparable to classical discontinuous Galerkin formulations in the respective discrete space norm on each mesh. Meanwhile, this technique allows us to solve on coarse meshes and adapt the solution to achieve a user-specified solution quality. 

\end{abstract}

\begin{keyword}
stabilized finite elements \sep residual minimization \sep inf-sup stability \sep advection-diffusion-reaction \sep adaptive mesh refinement
\end{keyword}

\end{frontmatter}

\section{Introduction}

Advection-diffusion-reaction problems arise in a wide range of phenomena relevant to many areas of applied physics and engineering, such as flow in porous media (e.g., reservoir engineering~\cite{ ewing2001, Calo2014} and groundwater flow~\cite{ Calo2014, ern2009}) and drug delivery~\cite{ Hossain2012, Calo2008, Bazilevs2007}. These processes generally involve heterogeneous and highly anisotropic diffusion tensors, representing varying material properties (e.g., permeability, porosity) in the domain~\cite{ Calo2016, galvis2018, Calo2011}. Thus, the accuracy and stability of the numerical approximation have been the focus of intense research for several decades. Moreover, in advection-dominated regimes, this problem behaves as a hyperbolic partial differential equation (PDE), implying that an inaccurate numerical approximation could produce non-physical oscillatory discrete solutions on coarse meshes. 

Traditionally, stabilization terms modify standard finite elements to improve the properties of the discrete solution. A plethora of conforming numerical schemes exist which add extra scaled residual terms~\cite{ brooks1982, hughes1989new, hughes2007}, and the residual-free bubble approximation (RFB)~\cite{ brezzi1998}. Despite these efforts, the stabilization process still represents a significant challenge for the scientific community due to its relevance in the outcome of numerical simulations. For instance, the discrete solution quality and solver performance strongly depend on the appropriate selection of the penalization (scaling) parameters.  Alternatively, non-conforming schemes build stability differently; for example, the discontinuous Galerkin (dG) class of methods, achieve stability by enforcing an element-by-element discretization and introducing a suitable choice of numerical traces~\cite{ reed1973, lesaint1974, johnson1986, cockburn2012, brezzi2004, ern2006}. Other conforming stabilized formulations are the minimal residual methods, such as the Least Squares Finite Element Method (LS-FEM)~\cite{ bochev2009} or the Discontinuous Petrov-Galerkin (DPG) method~\cite{ demkowicz2010, demkowicz2011, DemGopBOOK-CH2014, Calo2014dPG, Demkowicz2012, Niemi2011a, Niemi2013, Niemi2011}, which minimize the discrete residual with respect to an artificial energy norm and, thus, achieve the sought stability. 

Recently,~\cite{calo2019} introduced, analyzed, and numerically tested a stabilized finite element formulation in abstract form for any linear partial differential equation system. The procedure minimizes the discrete residual of a conforming trial space in the dual norm of a suitable dG test space in which the trial space is a proper subspace. Hence, this formulation endows the discrete, continuous solution with the stability properties of the dG formulation that defines the dual norm. The method expresses the constrained residual minimization as a saddle-point problem that delivers a conforming approximation and an on-the-fly error estimate that drives automatic mesh refinement. Besides, since the formulation builds on the non-conformity of the underlying dG method, this attribute, in turn, allows us to consider strong norms for the test space when the trial space has high regularity. Measuring the error in stronger norms is an outstanding feature of this technology when compared against alternative ones such as the LS-FEM and DPG methods. This non-conforming residual minimization successfully tackled the analysis of compaction banding in geomaterials~\cite{ cier2020}, the constraint enforcement in advection-dominated problems~\cite{ cier2020_2}, and incompressible flows~\cite{ los2020, kyburg2020}, as well as goal-oriented adaptive mesh refinements~\cite{ rojas2020}.

In this paper, we apply this stabilized finite element method to advection-diffusion-reaction problems and show the impact of the norm selection on the quality of the resulting solutions. We consider different possible scenarios for this class of elliptic problems, such as advection-dominated diffusion, and heterogeneous, with locally vanishing or highly-anisotropic diffusivities. Section~\ref{sec:problem} describes the advection-diffusion-reaction problem and a discontinuous Galerkin formulation that allows us to solve the problem via residual minimization. Section~\ref{sec:method} presents the method and the key ingredients for its formulation, and applies it to the model problem. Finally, Section~\ref{sec:examples} discusses several numerical examples that show the efficiency and robustness of this formulation.

\section{The advection-diffusion-reaction problem}\label{sec:problem}

\sloppypar
Let ${\Omega \subset \R^d}$, with dimension ${d = 2,3}$, be an open and bounded Lipschitz domain with boundary ${\Gamma :=\partial \Omega}$, and outward unit normal vector $\Vn$.  Using the standard notation of Hilbert and Banach spaces, let ${K \in \left[L^\infty(\Omega)\right]^{d,d}}$ be a diffusion tensor, to be symmetric and positive definite in $\Omega$. Let ${\Vb \in \left[L^\infty(\Omega)\right]^d}$ denote a divergence-free (almost everywhere) advection coefficient, and ${\sigma \in L^\infty(\Omega)}$ be a reactive coefficient. We write the advection-diffusion-reaction problem as follows:
\begin{equation}\label{eq:scalar}
  \left\{
    \begin{array}{l}
      \text{Find } u \text{ such that:} \smallskip \\
      \begin{aligned}
        -\div \left( K \nabla u \right) - \Vb \cdot \nabla u + \sigma \, u &= f,&& \text{ in } \Omega, \smallskip\\
        u &= 0, && \text{ on } \Gamma,
      \end{aligned}
    \end{array}
  \right.
\end{equation}
where ${f \in L^2(\Omega)}$ denotes an spatial source.

In the present work, we focus on two different types of advection-diffusion-reaction problems:
\begin{itemize}
\item Advection-dominated problems, that is, problems where ${{0<\|K\|_{\infty}, {\|\sigma\|_{\infty}}<<\|\Vb\|_{\infty}}}$. This kind of problems lead unstable solutions when using the standard FEM on coarse meshes.
\item Highly heterogeneous and anisotropic diffusion problems, that is, problems where the diffusion takes locally small values, leading to advection-dominated regimes.
\end{itemize}
These two scenarios lead to sharp inner and boundary layers, which are difficult to capture with standard FEM formulation as they induce spurious oscillations (see~\cite{ codina1998, hughes2018}).

\subsection{Continuous weak variational formulation}

The weak formulation of~\eqref{eq:scalar} reads:
\begin{equation}\label{eq:vf_cont}
  \left\{
    \begin{array}{l}
      \text{Find } u \in H^1_0(\Omega), \text{ such that:} \smallskip \\
      b(u,v)=\ell(v), \quad \forall v \in H^1_0(\Omega),
    \end{array}
  \right.
\end{equation}
with bilinear form $b(u,v)=(K\nabla u, \nabla v)_{0, \Omega} + (\Vb \cdot \nabla u, v)_{0, \Omega}+ (\sigma u, v)_{0, \Omega}$, and linear form $\ell(v)=(f, v)_{0, \Omega}$, where ${(\cdot, \cdot)_{0, \Omega}}$ denotes the $L^2$-scalar product in $\Omega$. In what follows, we assume that there is a real number ${\sigma_0 > 0}$ such that
\begin{align}
  \sigma - \frac{1}{2}\nabla \cdot \Vb \geq \sigma_0 \quad \textrm{a.e. in } \Omega.
\end{align}
Furthermore, we assume that the smallest eigenvalue of $K$ is bounded from below by a positive constant $K_0$. Then, owing to the Lax-Milgram Lemma, problem~\eqref{eq:vf_cont} is well-posed (see e.g.,~\cite{  ern2009}).

\subsection{Discrete setting}\label{ss:setting}

Let $\{\Par_h\}$ be a family of simplicial meshes of $\Omega$ and, for simplicity, we assume that any mesh exactly represents $ \Omega$ in $\Par_h$, that is, $\Omega$ is a polygon or a polyhedron. Let $T$ be a generic element in $\Par_h$, and denote by $\partial T$ its boundary, by $h_T$ its diameter, and by $\Vn_T$ its outward unit normal. We set ${h = \max_{T \in \Par_h} h_T}$. We define the classical dG approximation space as
\begin{align}
  V_h :=\{v_h \in L^2(\Omega) \ | \  \forall T \in \Par_h, v_h|_T \in \mathds{P}^p(T) \},
\end{align}
where $\mathds{P}^p(T)$ denotes the space of polynomial functions with degree smaller or equal than $p$. We set the extended space ${V_{h,\#} =H^2(\Par_h)+V_h}$ for convenience. We say that $F$ is an ``interior face'' if there exist two elements $\{T^-(F),\,{T^+(F)\} \in \Par_h}$, such that ${T^-(F) \cap T^+(F) = F}$ and $F$ has nonzero measure. We collect all faces $F$ of $\Par_h$ into the set ${\Sk_h= \bigcup_{T\in\Par_h} F}$. We denote by $\Sk^\partial_h$ the boundary skeleton ${\Sk^\partial_h = \Sk_h \cap \Gamma}$, and by $\Sk^0_h$ the internal skeleton ${\Sk^0_h=\Sk_h \backslash \Gamma}$.  Over each $F\in \Sk_h$, we set $\Vn_F$ as a predefined unit normal, oriented from $T^-(F)$ to $T^+(F)$, being coincident with $\Vn$ when ${F \in \Sk_h^\partial}$, and $h_F$ as the diameter of the face $F$. On interior faces, any function ${v_h \in V_h}$ is two-valued, with values $v_h^+$ and $v_h^-$, defined with respect to the predefined normal $\Vn_F$. Thus, the jump $\llbracket v_h \rrbracket_F$ and the weighted average $\llbrace v_h \rrbrace_\omega$ functions are defined as:
\begin{align*}
\llbracket v_h \rrbracket_F := v_h^+ - v_h^-, && \llbrace v_h \rrbrace_\omega:= \omega^- v_h^- + \omega^+v_h^+,
\end{align*}
where the weights satisfy ${\omega^- + \omega^+ =1}$, with ${\omega^-,\omega^+\geq 0}$. In particular, when considering heterogeneous tensorial diffusivities, we choose the weights accounting the diffusivity structure:
\begin{align}\label{eq:omega}
\omega^- = \frac{\delta_{Kn}^+}{\delta_{Kn}^+ + \delta_{Kn}^-}, \quad \quad \omega^+ = \frac{\delta_{Kn}^-}{\delta_{Kn}^+ + \delta_{Kn}^-},
\end{align}
with $\delta^\mp_{Kn}=\Vn_F\cdot K^\mp \Vn_F$ if $F \in \Sk_h^0$, and $\delta_{Kn}=\Vn_F\cdot K \Vn_F$ if $F \in \Sk_h^\partial$. When $K$ is a continuous tensor (homogeneous diffusion), the weights reduce to ${\omega^-=\omega^+=1/2}$. Finally, on a boundary face ${F \in \Sk^\partial_h}$, we set ${\llbracket v_h \rrbracket_F=\llbrace v_h \rrbrace_F=v_h|_F}$. We omit the subscript $F$ in the jump and weighted average functions when there is no ambiguity.

\subsection{Discontinuous Galerkin variational formulation}\label{section:dg_form}

In this section, we briefly discuss a stable discontinuous Galerkin (dG) formulation for problem~\eqref{eq:scalar}. It combines the Symmetry Weighted Interior Penalty (SWIP) scheme~\cite{ di2008, ern2009} that handles general diffusivities, combined with the Upwinding (UPW) method~\cite{ brezzi2004, di2011mathematical} that handles the advection-reaction contribution. 

Considering the discrete setting described in Section~\ref{ss:setting}, the dG formulation of problem~\eqref{eq:vf_cont} reads:
\begin{equation}\label{eq:gen_vf}
  \left\{
    \begin{array}{l}
      \text{Find } u_h \in V_h, \text{ such that:} \smallskip \\
      \begin{aligned}
        b_h(u^{\text{dG}}_h,v_h)&:=b_h^\textrm{diff}(u^{\text{dG}}_h,v_h) + b_h^\textrm{adv}(u^{\text{dG}}_h,v_h)= \ell_h(v_h), && \forall v_h \in V_h,
      \end{aligned}
    \end{array}
  \right.
\end{equation}
%
where the bilinear and linear forms are
\begin{align*}
  b_h^\textrm{diff}(u_h,v_h)&:=  \displaystyle \sum_{T \in \Par_h} (K \nabla u_h \, , \, \nabla v_h)_{0,T} \\
                            & + \displaystyle \sum_{F \in \Sk_h} \left[ \displaystyle \, - \left(  \llbracket u_h \rrbracket \, , \, \Vn_F\cdot \llbrace  K \nabla v_h \rrbrace_{\omega} \right)_{0,F} - \displaystyle \left( \Vn_F\cdot \llbrace  K \nabla u_h \rrbrace_{\omega} \, , \,  \llbracket v_h \rrbracket \right)_{0,F}  + \gamma_F \left( \llbracket u_h \rrbracket \, , \, \llbracket v_h \rrbracket \right)_{0,F} \right] , \smallskip\\
  b_h^\textrm{adv}(w,v) & :=  \displaystyle \sum_{T\in\Par_h}(\, \Vb \cdot \nabla  u_h + \sigma \,  u_h\, , \, v_h)_{0,T} + \displaystyle \sum_{F\in\Sk^\partial_h} \left( (\, \Vb \cdot \Vn_F)^\ominus u_h, v_h\right)_{0,F} \\
                            &+ \sum_{F \in \Sk^0_h} \left[ \dfrac{1}{2} \left(\, \left| \, \Vb \cdot \Vn_F\right| \, \llbracket u_h \rrbracket \, , \, \llbracket v_h \rrbracket \right)_{0,F} - \left( \, \Vb \cdot \Vn_F \, \llbracket u_h \rrbracket \, , \, \llbrace v_h \rrbrace \right)_{0,F}  \right], \smallskip
\end{align*}
and
\begin{align*}
\ell_h(v_h) &:= \sum_{T\in\Par_h}(f, v_h)_{0,T} 
\end{align*}
In the above, $(\cdot)^{\ominus}$ denotes the negative part of $x$ (i.e., ${x^\ominus := \frac{1}{2}\left(|x|-x\right)}$ for any real number $x$). For problems with diffusion, there exist many suitable choices for the penalty parameter $\gamma_F$ (e.g.,~\cite{ shahbazi2005, epshteyn2007} analyzed the penalty parameters and their dependence on the polynomial order of approximation;~\cite{ ern2009} defined and analyzed the impact of $\gamma_K$ as the harmonic average of the ``normal" permeabilities; ~\cite{ hartmann2008}, introduced a mesh-dependent penalty parameter). In this work, following~\cite{ern2009, shahbazi2005, bastian2012}, we set the penalty parameter $\gamma_F$ as ${\gamma_F = \displaystyle \eta \,{\gamma_K}}$. Here, $\eta>0$ represents an element-wise parameter defined as:
  \begin{align}
    \forall F \in \Sk^0_h, && \eta &= \frac{1}{2}\frac{(p+1)(p+d)}{d}\left(\frac{\Ak(\partial T^+)}{\Vk(T^+)}+\frac{\Ak(\partial T^-)}{\Vk(T^-)}\right), \\
    \forall F \in \Sk^\partial_h, && \eta &= \frac{(p+1)(p+d)}{d}\frac{\Ak(\partial T)}{\Vk(T)},
  \end{align}
  where $p$ is the polynomial order of the test space $V_h$, $d$ is the dimension, and $\Ak$ and $\Vk$ denote area and volume, respectively, for $d=3$, and length and area, respectively, for $d=2$. We define ${\gamma_K}$ as follows:
  \begin{align}
    \forall F \in \Sk^0_h, && \gamma_K &= (\omega^-)^2\delta^-_{Kn} + (\omega^+)^2\delta^+_{Kn}, \\
    \forall F \in \Sk^\partial_h, && \gamma_K &= \delta_{Kn},
  \end{align}
  thus, recalling weights definition~\eqref{eq:omega}, we derive that
  \begin{align}
    \forall F \in \Sk_h^0, &&\gamma_K &= \frac{\delta_{Kn}^+ \delta_{Kn}^-}{\delta_{Kn}^+ + \delta_{Kn}^-}.
  \end{align}

When we consider scalar diffusivities, that is, $K = \varkappa I$ for some scalar function ${\varkappa:\Omega \rightarrow \R}$, we recover the symmetric interior penalty (IP or SIP) method~\cite{ wheeler1978, arnold1982, arnold2002unified} given that the penalty parameters reduces to $\gamma_F = \eta \varkappa$.

\begin{rmrk}
  The dG formulation allows us to weakly impose non-homogeneous Dirichlet boundary conditions through modifying the right-hand side of~\eqref{eq:gen_vf} as follows: if we look for a solution of problem~\eqref{eq:scalar} satisfying $u=g_D$ on $\Gamma$, being ${g_D \in H^{1/2}(\Gamma)}$ a boundary source, then we rewrite the linear form $\ell_h(v_h)$ as:
  \begin{align*}
    \displaystyle \ell_h(v_h) &:= \displaystyle \sum_{T\in\Par_h}(f, v_h)_{0,T} +  \sum_{F \in \Sk_h^\partial}  \left[ \, - \left( g_D,  \Vn_F\cdot K \nabla v_h \right)_{0,F} + \gamma_F \left( g_D, v_h \right)_{0,F}   + \left( ( \, \Vb \cdot \Vn_F)^\ominus g_D,  v_h\right)_{0,F} \right]. \smallskip
  \end{align*}
\end{rmrk}

\subsection{Discrete norms, well-posedness and a priori error estimate}\label{ss:discrete_norms}

We provide the discrete space $V_h$ with the following norm: 
\begin{align}\label{eq:norm}
  \|w\|^2_{V_h} := \|w\|^2_{\text{adv}} + \|w\|^2_{\text{diff}}
\end{align}
with
\begin{align*}
  \|w\|^2_{\text{adv}} &:= \displaystyle \|w\|^2_{0, \Omega} + \frac{1}{2}\| \, | \, \Vb\cdot \Vn \, |^{\frac{1}{2}}w\|^2_{0, \Gamma} + \displaystyle \frac{1}{2}\sum_{ F \in \mathscr{S}^0_h}\left(| \, \Vb\cdot \Vn_F | \, \llbracket w \rrbracket, \llbracket w \rrbracket \right) _{0,F} + \displaystyle \sum_{T \in \mathscr{P}_h} h_T \| \, \Vb\cdot\nabla{w} \, \|^2_{0, T}, \\
  \|w\|^2_{\text{diff}} &:= \displaystyle \| \, \kappa \nabla w \, \|^2_{0, \Omega} +  \sum_{F \in \mathscr{S}_h} \left( \gamma_{F}\llbracket w \rrbracket, \llbracket w \rrbracket \right) _{0,F},
\end{align*}
where $\kappa$ denotes the (unique) symmetric positive definite tensor-valued field such that $\kappa^2 = K$ a.e. in $\Omega$. Following~\cite{ ern2009}, we define:
\begin{equation}\label{eq:Vh_beta}
  |w|^2_{V_{h,\beta}} := \displaystyle \sum_{T \in \mathscr{P}_h} h_T \| \, \Vb\cdot\nabla{w} \, \|^2_{0, T},
\end{equation}
which represents the last component of the norm $\|w\|^2_{\text{adv}}$ that controls the advective derivative error for small diffusivities. Finally, we define the following extended norm $\|w\|_{V_h, \#}$:
\begin{equation}\label{eq:Vh_ext}
\|w\|_{V_h, \#} := \displaystyle \|w\|_{V_h} + \left( \sum_{T \in \Par_h}  \| w \|^2_{0,\partial T} \right)^{\frac{1}{2}} + \displaystyle \left( \sum_{T \in \Par_h} h_T \|\kappa \nabla w \|^2_{0,\partial T} \right)^{\frac{1}{2}}.
\end{equation}	

In the remainder, the symbol $\lesssim$ indicates an inequality involving a positive constant $C$ independent of the mesh and diffusivity. Considering the norms above defined, the following theorem holds true (see~\cite[\S3 \& \S4]{ ern2009}):
\begin{thrm}[Well-posedness and a priori error estimate of the dG formulation]\label{lmm:well-posedness}
  The following propositions hold true:
  \begin{enumerate}[label=(\alph*)]
    
  \item Inf-sup stability: There exists a constant ${C_{\emph{sta}}= C\Delta_K^{-1}}$, with ${C>0}$, uniform with respect to the mesh size, such that:
    \begin{align*}
      \sup_{v_h \in V_h \backslash \{0\}}\frac{b_h(w_h, v_h)}{\|v_h\|_{V_h}} \geq C_{\emph{sta}}\|w_h\|_{V_h}, \quad \forall w_h \in V_h,
    \end{align*}
    where $\forall \, T \in \Par_h$, ${\Delta_K =\max_{T\in\Par_h}\Delta_{K, T}}$, and
    \begin{align*}
      \quad \Delta_{K, T} = \left\{
      \begin{array}{l l}
        \smallskip 1 & \text{if} \ \  \|\Vb\|_{[L^\infty(T)]^d} \gtrsim \dfrac{\lambda_{M,T}}{h_T},\\
        \dfrac{\lambda_{M,T}}{\lambda_{m,T}} & \text{otherwise},
      \end{array}
                                               \right.
    \end{align*}
    with ${\lambda_{M,T}}$ and ${\lambda_{m,T}}$ as the maximum and minimum eigenvalues of $K|_T$, respectively. 

  \item Boundedness: There exists a mesh-independent constant $C_{\emph{bnd}} < \infty$, such that:
    \begin{align*}
      b_h(z, v_h) \leq C_{\emph{bnd}} \|z\|_{V_h,\#}\|v_h\|_{V_h}, \quad \forall (z,v_h) \in V_{h,\#} \times V_h.
    \end{align*}
		
  \item Consistency: Let $u$ be the solution of~\eqref{eq:vf_cont}. If ${u \in H^1_0(\Omega) \cap H^2(\Par_h)}$, then
    \begin{align*}
      b_h(u, v_h)=\ell_h(v_h), \quad \forall v_h \in V_h.
    \end{align*}
  \end{enumerate}
  Henceforth, $u^{\emph{dG}}$, solution of problem~\eqref{eq:gen_vf}, is unique (cf.~\cite{ern2009,di2011mathematical}). Moreover, the following a priori error estimate is satisfied
  \begin{align}\label{eq:apriori}
\inf_{y_h \in V_h} \|u - y_h\|_{V_h} \leq    \|u-u^{\emph{dG}}_h\|_{V_h} \leq \left(1 + \frac{C_\emph{bnd}}{C_\emph{sta}}\right) \inf_{v_h \in V_h} \|u - v_h \|_{V_h, \#}.
  \end{align}
  
\end{thrm}

We relate the convergence rates of both left and right-hand sides of the error estimate~\eqref{eq:apriori} through the following definition (cf.~\cite{ di2011mathematical}, \S1.4.4):
  \begin{dfntn}[Optimality, quasi-optimality and suboptimality of the error estimate]\label{def:opt}
    The error estimate~\eqref{eq:apriori} is
    \begin{enumerate}
    \item \emph{optimal} if $\| \cdot \|_{V_h} \simeq \| \cdot \|_{V_h, \#}$, 
    \item \emph{quasi-optimal} if the two norms are different, but the lower and upper bounds in~\eqref{eq:apriori} converge, for smooth $u$, at the same convergence rate as $h \rightarrow 0$,
    \item \emph{suboptimal} if the upper bound converges at a slower rate than the other bound.
    \end{enumerate}
  \end{dfntn}

\section{Residual minimization based on discontinuous Galerkin formulations}\label{sec:method}

\subsection{Method description}

We now describe the stabilized finite element formulation via residual minimization on dual discontinuous Galerkin (dG) norms devised in~\cite{ calo2019} for the advection-diffusion-reaction problem context.  In~\eqref{eq:gen_vf}, $u^{\text{dG}}_h$ discretely approximates the solution of~\eqref{eq:vf_cont} belonging to a \textit{discontinuous} discrete space. In~\cite{ calo2019}, the formulation seeks for an approximation of $u$ in a discrete space that may possess additional properties, for instance, continuity and, possibly, higher smoothness. Rather than solving problem~\eqref{eq:gen_vf}, the dG-based residual minimization method implies:
\begin{enumerate}[label=(\alph*)]
\item Consider the subspace $U_h = V_h \cap U \subset V_h$ (i.e., standard finite element functions).
\item Obtain ${u_h \in U_h}$ by minimizing the residual:
  \begin{equation}\label{eq:min_prob}
    \left\{
      \begin{array}{l}
        \text{Find } u_h \in U_h \subset V_h,  \text{ such that:} \smallskip \\
        u_h = \displaystyle \argmin_{z_h \in U_h} \dfrac{1}{2}\|\ell_h- B_h \, z_h\|^2_{V_h^\ast} = \displaystyle \argmin_{z_h \in U_h} \dfrac{1}{2}\|R^{-1}_{V_h}(\ell_h- B_h z_h)\|^2_{V_h},
      \end{array}
    \right.
  \end{equation}
  where $B_h : V_h, \# \rightarrow V_h^{\ast}$ is: 
  \begin{align}
    \langle B_h w_h, v_h \rangle_{V_h^{\ast} \times V_h} := b_h(w_h,v_h),
  \end{align}
  and $R^{-1}_{V_h}$ denotes the inverse of the Riesz map:
  \begin{align}\label{eq:riesz}
    R_{V_h} :\ & V_h \rightarrow V_h^\ast \smallskip\\
               & v_h \rightarrow \langle R_{V_h} y_h, v_h\rangle_{V_h^\ast \times V_h} := (y_h,v_h)_{V_h},
  \end{align}
  with $(\cdot, \cdot)_{V_h}$ denoting the inner product that induces the discrete norm ${\|\cdot\|_{V_h}}$ (i.e., ${\|\cdot\|_{V_h}=(\cdot, \cdot)^{1/2}_{V_h}}$).
\end{enumerate}

The second equality in~\eqref{eq:min_prob} follows as the Riesz map~\eqref{eq:riesz} is an isometric isomorphism, therefore $\| \cdot \|_{V_h^\ast}$ is equivalent to $\| R^{-1}_{V_h} (\cdot) \|_{V_h}$.  Thus, problem~\eqref{eq:min_prob} is equivalent to the following saddle-point problem (see~\cite{ cohen2012}):
\begin{equation}\label{eq:saddle-point}
  \left\{
    \begin{array}{l}
      \text{Find } (\varepsilon_h, u_h) \in V_h \times U_h,  \text{ such that:} \smallskip \\
      \begin{array}{rcll}
        (\varepsilon_h,v_h)_{V_h} + b_h(u_h,v_h) &=& \ell_h(v_h),& \quad \forall v_h \in V_h, \smallskip\\
        b_h(z_h, \varepsilon_h)& =&0,&\quad \forall z_h \in U_h,
      \end{array}
    \end{array}
  \right.
\end{equation}
where $\varepsilon_h$ is a residual representative in $V_h$. Indeed, the first identity in~\eqref{eq:saddle-point} implies that:
\begin{align}\label{eq:epsilon_h}
  \varepsilon_h = R^{-1}_{V_h} (\ell_h-B_h u_h)=R^{-1}_{V_h} B_h(u_h^{\textrm{dG}} - u_h),
\end{align}
where the second identity in~\eqref{eq:epsilon_h} comes from~\eqref{eq:gen_vf}.

The saddle-point formulation has several desirable properties for numerical approximations. Firstly, the matrix associated with the inner product $(\varepsilon_h,v_h)_{V_h}$ in~\eqref{eq:saddle-point} is always symmetric and positive-definite, independently of the nature of the chosen dG formulation; thus, several well-known iterative solvers are effective on the resulting saddle-point problem. Moreover, the discrete approximation $u_h \in U_h$ inherits the discrete stability of the dG formulation. Finally, $\varepsilon_h \in V_h$ is a robust error representation, which under an adequate saturation assumption becomes a reliable representative. Below, we summarize the last properties (see~\cite{ calo2019} for details):  
\begin{thrm}[Well posedness and a priori error bound estimates for the saddle-point problem]\label{thrm:mixed_wp}
  The solution $(\varepsilon_h, u_h) \in V_h \times U_h$ of the saddle-point problem~\eqref{eq:saddle-point} is unique and the following a priori bound applies:
  \begin{align}
    \|\varepsilon_h\| \leq \|\ell_h\|_{V_h^*} \quad \text{and} \quad \|u_h\|_{V_h} \leq \frac{1}{C_{\emph{sta}}}\|\ell_h\|_{V_h^*},
  \end{align}
  and the following error estimate holds true:
  \begin{equation}\label{eq:error-est}
    \|u-u_h\|_{V_h} \leq \left(1 + \dfrac{C_\emph{bnd}}{C_\emph{sta}} \right) \Delta_K  \inf_{v_h \in U_h}\|u-v_h\|_{V_h,\#},
  \end{equation}
  where ${u \in V_\#}$ represents the exact solution to the continuous problem~\eqref{eq:scalar}. 

\end{thrm}

\begin{prpstn}[Efficiency of the residual representative]
  Under the same hypotheses of Theorem~\ref{thrm:mixed_wp}, the following holds:
  \begin{equation}\label{eq:efficient}
    \|\varepsilon_h\|_{V_h}\leq C_{\emph{bnd}} \, \|u - u_h \|_{V_h, \#}.
  \end{equation}
\end{prpstn}

\begin{ssmptn}[Saturation]\label{as:saturation}
  Let ${u_h\in U_h}$ be the second component of the pair ${(\varepsilon_h,u_h)\in V_h\times U_h}$ solving the saddle-point problem~\eqref{eq:saddle-point}. Let ${u^{\emph{dG}}_h\in V_h}$ be the unique solution to~\eqref{eq:gen_vf}. There exists a real number $\delta\in [0,1)$, uniform with respect to the mesh size, such that ${\|u-u^{\emph{dG}}_h\|_{V_h} \le \delta \|u-u_h\|_{V_h}}$.
\end{ssmptn}

\sloppy \begin{prpstn}[Reliability of the residual representative] \label{th:err_rep}
  Let ${u_h\in U_h}$ be the second component of ${(\varepsilon_h,u_h)\in V_h\times U_h}$ solving the saddle-point problem~\eqref{eq:saddle-point}. Let ${u^{\text{dG}}_h \in V_h}$ be the unique solution to~\eqref{eq:gen_vf}. Then the following holds true:
  \begin{equation} \label{eq:bnd_uh-thetah}
    \|u_h-u^{\emph{dG}}_h\|_{V_h} \le \dfrac{1}{C_{\emph{sta}}} \|\varepsilon_h\|_{V_h}.
  \end{equation} 
  Moreover, if the saturation Assumption~\ref{as:saturation} is satisfied, then the following a posteriori error estimate holds:
  \begin{equation}\label{eq:a_posteriori}
    \|u-u_h\|_{V_h} \le \dfrac{1}{(1-\delta)C_{\emph{sta}}} \|\varepsilon_h\|_{V_h}.
  \end{equation} 
\end{prpstn}

We can now state the requirements for having an efficient and reliable residual representative. From~\eqref{eq:efficient} and~\eqref{eq:a_posteriori}, we have:
  \begin{equation}\label{eq:rates}
    \|u-u_h\|_{V_h} \lesssim \|\varepsilon_h\|_{V_h} \lesssim \|u - u_h \|_{V_h, \#}.
  \end{equation}
Thus, to ensure the usefulness of the residual representative, we need at least quasi-optimality, which means that the left-hand side should decay at the same rate as the right-hand side of~\eqref{eq:rates} (see Definition~\ref{def:opt}).

\subsection{Convergence rates}

In the context of the above framework, we can recover some insightful results related to the convergence rates. Similarly as done in~\cite[Appendix B]{calo2019} (see also~\cite{ karakashian2003, burman2007, ern2017}), we can prove:
\begin{equation}\label{eq:beta-1}
\inf_{v_h \in U_h}\|u-v_h\|_{V_h,\#} \, \lesssim \, \inf_{v_h \in V_h} \|u-v_h\|_{V_h,\#}.
\end{equation}
This proves, in particular, that the residual minimization method delivers a discrete solution with the same quality as the one delivered by the dG formulation, as consequence of~\eqref{eq:apriori}~and~\eqref{eq:error-est}. Thus, if the solution is regular enough, it follows:
\begin{equation}\label{eq:beta-3}
\inf_{v_h \in V_h} \|u-v_h\|_{V_h,\#} \lesssim h^{p}.
\end{equation}
Additionally, from~\eqref{eq:beta-3}  we can deduce a bound for the advective component ${|\cdot|_{V_h,\beta}}$, defined in~\eqref{eq:Vh_beta}. Indeed, we can easily infer that:
\begin{equation}\label{eq:beta-2}
  \| \, \Vb\cdot\nabla(u-u_h) \, \|^2_{0, T} \, \lesssim \, \| \, \kappa \nabla (u-u_h) \, \|^2_{0, \Omega} \, \lesssim \, \|u-u_h\|^2_{V_h,\#}, 
\end{equation}
where $u$ and $u_h$ are defined as in Theorem~\ref{thrm:mixed_wp}. 
Combining~\eqref{eq:beta-1},~\eqref{eq:beta-2} and~\eqref{eq:beta-3}, multiplying both sides by the mesh size, and taking the square root, we can recover the optimal convergence rate for ${|u-u_h|_{V_{h,\beta}}}$, which reads:
\begin{equation*}
  |u-u_h|_{V_{h,\beta}}\, = \left(\displaystyle \sum_{T \in \mathscr{P}_h} h_T \| \, \Vb\cdot\nabla(u-u_h) \, \|_{0, T}\right)^{1/2} \, \lesssim \, h^{1/2} \|u-u_h\|_{V_h,\#} \, \lesssim \, h^{p+1/2}.
\end{equation*}
Similarly, we can deduce a bound for the $L^2$-norm  error in the following way:
\begin{equation}\label{eq:l2}
  \|  u-u_h \|_{0, \Omega} \, \lesssim \, \|u-u_h\|_{V_h,\#} \, \lesssim \, \inf_{v_h \in V_h} \|u-v_h\|_{V_h,\#} \lesssim h^p,
\end{equation}
which is suboptimal. Indeed, we show through numerical examples that the method could deliver solutions with suboptimal $L^2$ convergence  rates. Nevertheless, we show that the optimal convergence can be recovered by performing the energy-norm based adaptive mesh refinement strategy.

\section{Numerical examples}\label{sec:examples}

In this section, we discuss some implementation aspects and describe the test cases that demonstrate the performance of the method under a wide range of challenging scenarios.

\subsection{Implementation aspects}

We use FEniCS~\cite{ alnaes2015fenics} to perform all the numerical simulations. We show convergence plots of the error measured in the $L^2$ and $V_h$ norms versus the number of degrees of freedom (DOFs) (i.e., $\dim(U_h) + \dim(V_h)$ for~\eqref{eq:saddle-point}). As Section~\ref{section:dg_form} describes, we use the Symmetric Weighted Interior Penalty (SWIP) method for the diffusive part of the bilinear form, which extends the classical Symmetric Interior Penalty (SIP) formulation. In our experience, extensive numerical testing shows that other dG formulations have similar computational cost and convergence rates; thus, for the sake of brevity, we detail the performance of the method uniquely using the SWIP formulation.

\subsubsection{Adaptive mesh refinement}

We use $\varepsilon_h\in V_h$ to estimate the error and to drive the adaptive mesh refinement process~\cite{ calo2019}. We follow a standard adaptive procedure, which considers an iterative loop where each level of refinement we perform the following four analysis modules:
$$
\text{ SOLVE } \rightarrow \text{ ESTIMATE } \rightarrow \text{ MARK } \rightarrow  \text{ REFINE. } 
$$
That is, given a mesh partition, we first solve the saddle-point problem~\eqref{eq:saddle-point}. Later, we use a localized version of the inner product~\eqref{eq:norm} (evaluated in each mesh cell $T$) as error indicator $E_T$
\begin{equation}
E_T ^2 := \| \varepsilon_h \|^2_{loc, T} + \frac12 |\varepsilon_h|^2_{loc,S} \, ,
\end{equation}
with
\begin{equation}
\begin{array}{rcl}
\| \varepsilon_h \|^2_{loc, T} &:=& \, \displaystyle \|\varepsilon_h\|_{0,T}^2 + \displaystyle \|\kappa \nabla \varepsilon_h \|_T^2 + \displaystyle h_{0,T} \|\ \Vb \cdot \nabla \varepsilon_h\|^2_{0,T} \displaystyle + \sum_{F \in \mathscr{S}^\partial_h} \left( {\frac12}  \, | \Vb \cdot \Vn |  + \gamma_F \right) \left( \varepsilon_h , \varepsilon_h \right) _{0,F}, \\
|\varepsilon_h|^2_{loc,S} &:=& \displaystyle \sum_{F \in \mathscr{S}^0_h} \left( {\frac12}  \, | \Vb \cdot \Vn |  + \gamma_F \right) \left( \llbracket \varepsilon_h \rrbracket, \llbracket \varepsilon_h \rrbracket \right) _{0,F}.
\end{array}
\end{equation}
We then mark elements for refinement using the D{\"o}rfler bulk-chasing criterion~\cite{ dorfler1996convergent} (i.e., select all elements for which the cumulative sum of the local values $E_T$ in a decreasing order remains below a user-defined fraction of the total estimated error $\|\varepsilon_h\|_{V_h}$). Herein, we set this fraction to be 0.5 for $d=2$ and 0.25 for $d=3$. Finally, we bisect each marked element~\cite{ bank1983some} to obtain the refined mesh to use in the next step.

\subsubsection{Iterative solver}

The matrix system that problem~\eqref{eq:saddle-point} induces has the following form
\begin{equation}
  \left[
    \begin{array}{cc}
      G & B \\
      B^T & 0
    \end{array}
  \right] 
  \left[
    \begin{array}{c}
      {\boldsymbol{\varepsilon}}\\
      {\bf u}
    \end{array}
  \right] =
  \left[
    \begin{array}{c}
      \bf{L} \\
      {\bf 0} 
    \end{array}
  \right].
\end{equation}
where the superindex $T$ denotes transpose. Following~\cite[\S5]{ calo2019}, we apply the iterative algorithm proposed in~\cite{ bank1989class}. Denoting by $\widehat{G}$ a preconditioner for the Gram matrix $G$, and by $\widehat{S}$ a preconditioner for the reduced Schur complement $B^T \widehat{G}^{-1} B$, the iterative scheme becomes
\begin{equation}\label{eq:iter_solver}
  \left[
    \begin{array}{c}
      {\boldsymbol \varepsilon}_{i+1}\\
      {\bf u}_{i+1}
    \end{array}
  \right]
  = 
  \left[
    \begin{array}{c}
      {\boldsymbol \varepsilon}_i\\
      {\bf u}_i
    \end{array}
  \right]
  + 
  \left[
    \begin{array}{cc}
      \widehat{G} & B \\
      B^T & \widehat{C}
    \end{array}
  \right]^{-1} 
  \left\{
    \left[
      \begin{array}{c}
        \bf{L} \\
        {\bf 0} 
      \end{array}
    \right]
    -
    \left[
      \begin{array}{cc}
        G & B\\
        B^T & 0
      \end{array}
    \right] 
    \left[
      \begin{array}{c}
        {\boldsymbol \varepsilon}_i\\
        {\bf u}_i
      \end{array}
    \right]
  \right\},
\end{equation}
with $\widehat{C} = B^T \widehat{G}^{-1} B - \widehat{S}$. Let ${{\bf r}_i= {\bf L} - G \, {\boldsymbol \varepsilon}_i - B\, {\bf u}_i}$ and ${{\bf s}_i = - B^T {\boldsymbol \varepsilon}_i}$. Also, let ${\boldsymbol \varepsilon}$ and ${\bf u }$ be the residuals at the outer iteration $i$. Then, the scheme requires the resolution of two interior problems for these increments:
\begin{equation}\label{eq:iter_solver_1}
  \boldsymbol \eta_{i+1}:={\bf u}_{i+1}-{\bf u}_i = \widehat{S}^{-1}\left(B^T\left(\widehat{G}^{-1} {\bf r}_i\right) - {\bf s}_i \right),
\end{equation}
and
\begin{equation}\label{eq:iter_solver_2}
\boldsymbol\delta_{i+1}:={\boldsymbol\varepsilon}_{i+1}-{\boldsymbol\varepsilon}_{i} = \widehat{G}^{-1}\left({\bf r}_i - B \, {\boldsymbol \eta}_{i+1} \right).
\end{equation}
We construct an accurate approximation of $G^{-1}$, as in practice, this delivers the best good computational performance; that is, low-quality approximations lead poor conditioning of the reduced Schur complement in~\eqref{eq:iter_solver_1}. \citet{ bank1989class} proposed a relaxed approximation for the matrix $G$; nevertheless, for this stiffer problem, we use a scheme with one outer iteration. That is, we approximate $G^{-1}$ with a sparse Cholesky factorization (e.g., using the module ``sksparse.cholmod'', see~\cite{ chen2008algorithm}). Moreover, we choose precondition~\eqref{eq:iter_solver_1} with an approximate Schur complement built as ${\widehat{S}=B^T(\text{diag}(G))^{-1}B}$, where ${\text{diag}(G)}$ is the main diagonal of $G$, and the inverse is approximated through the same procedure used for $G$. We use the LGMRES algorithm (e.g., from Scipy sparse linear algebra package) to solve. On the coarsest mesh, our initial guess is zero, whereas, on the next refinements, our guess is the solution of the previous level of refinement.

%
%
%

\begin{figure}[h!]
	\begin{center}
		\begin{subfigure}{0.325\textwidth}
			\centering
			\includegraphics[width=\textwidth]{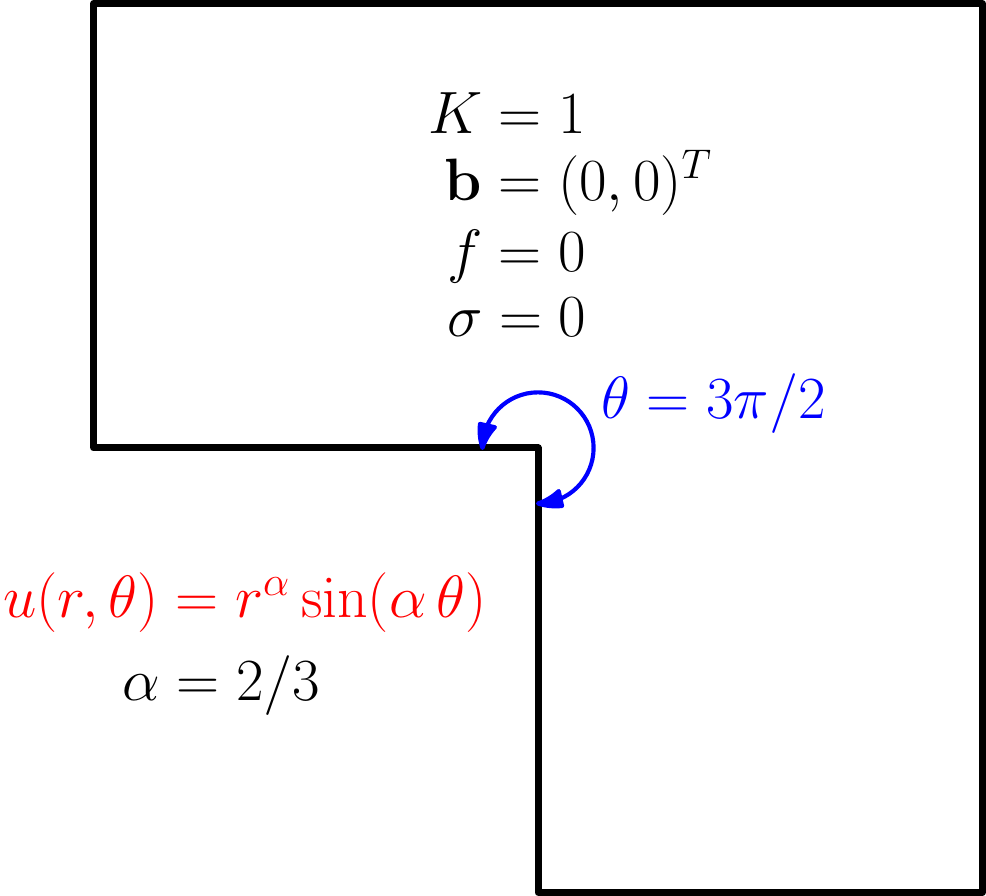} 
			\label{L-shape-23:skecth}
			\caption{Problem sketch}
		\end{subfigure}
	\hspace{0.15\textwidth}
		\begin{subfigure}{0.425\textwidth}
			\centering
			\includegraphics[width=\textwidth]{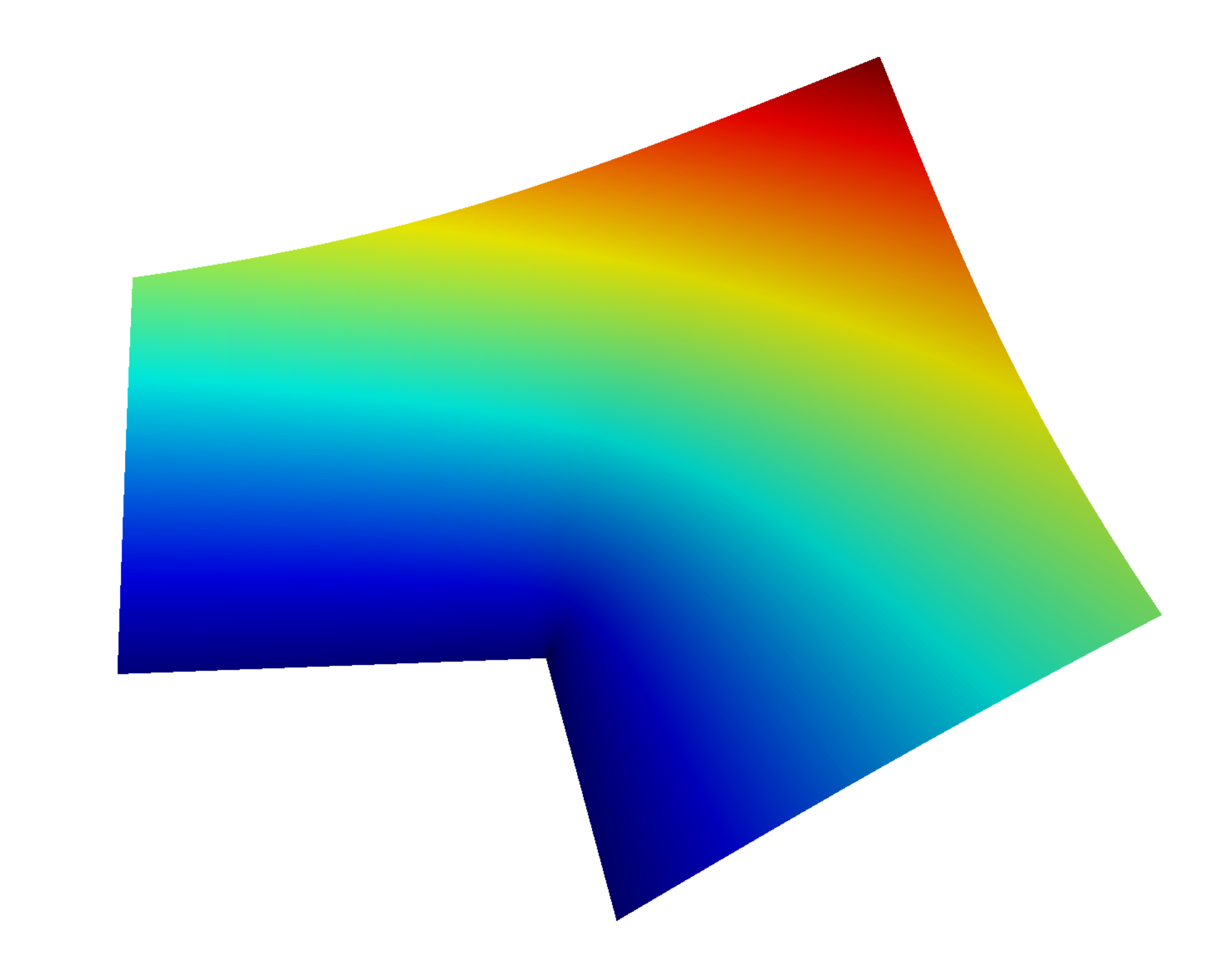} 
			\label{L-shape-23:3D}
			\vspace{-0.05\textwidth}
			\caption{3D view of the adaptive solution for $p=1$}
			\vspace{0.05\textwidth}
		\end{subfigure}
	\hspace{0.05\textwidth}
		\begin{subfigure}{0.6\textwidth}
			\centering
			\includegraphics[width=.8\textwidth]{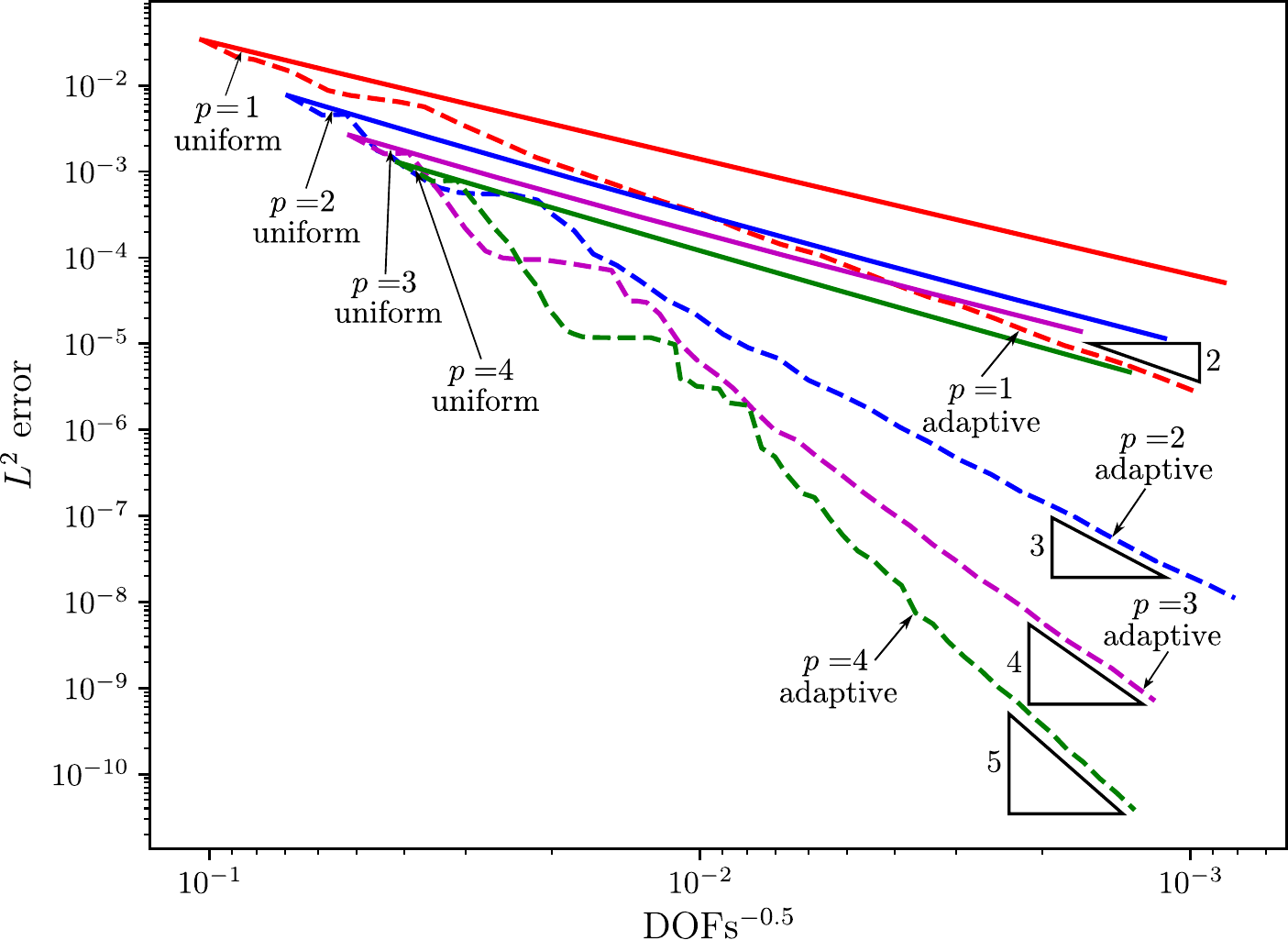} 
			\label{L-shape-23:L2}
			\caption{Convergence plot in $L^2$ error norm. $\alpha=2/3$}
		\end{subfigure}
	\hspace{0.05\textwidth}
		\begin{subfigure}{0.6\textwidth}
			\centering
			\includegraphics[width=.8\textwidth]{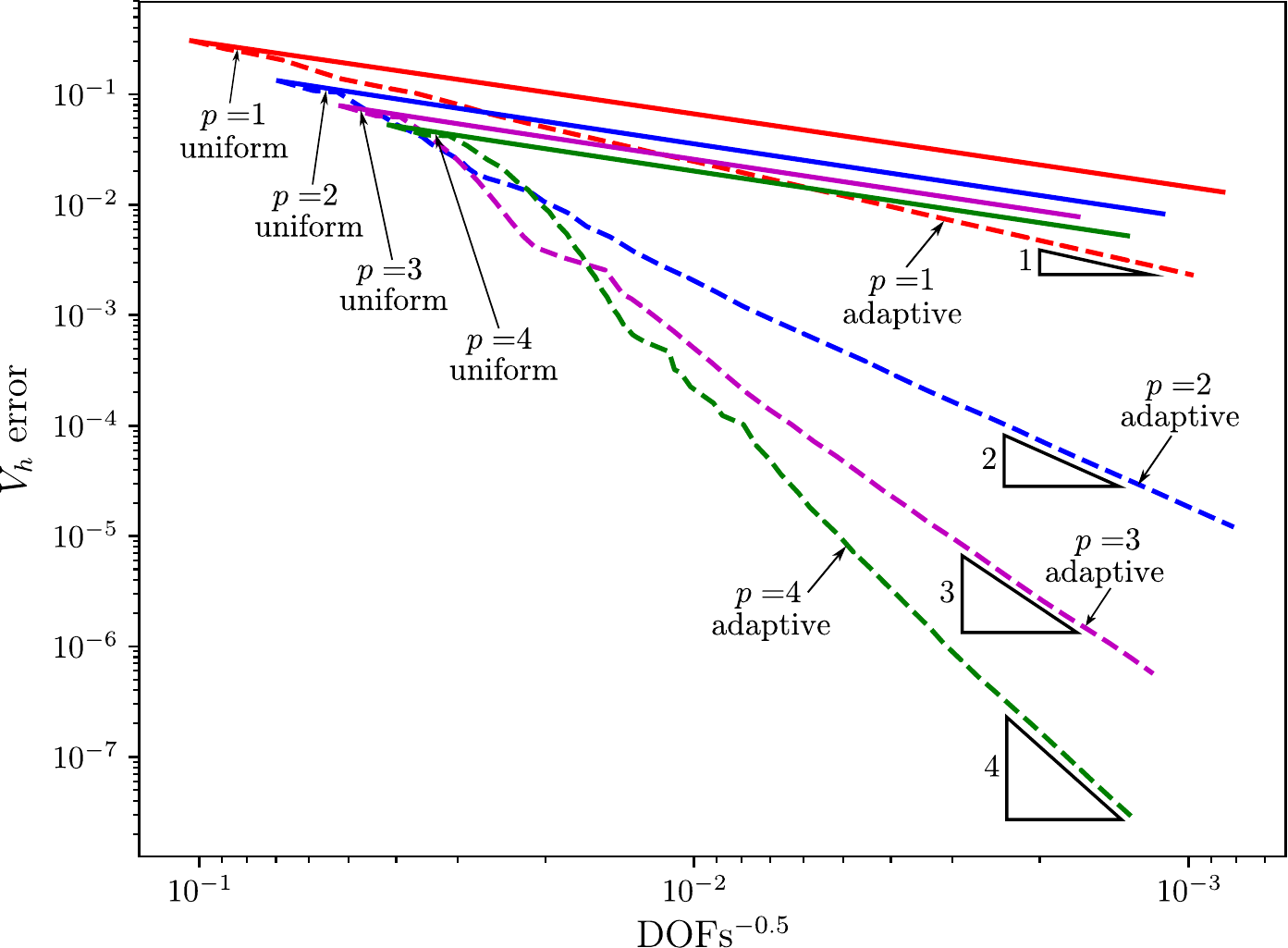} 
			\label{L-shape-23:Vh}
			\caption{Convergence plot in $V_h$ error norm. $\alpha=2/3$}
		\end{subfigure}
	\end{center}
	\caption{Convergence rates,   Poisson problem in L-shape domain. $\alpha=2/3$}
	\label{fig:L-shape-43}
\end{figure}

\subsection{Pure diffusion on L-shape domain}

As a first example, we consider the L-shape domain ${\Omega := (-1,1)^2 \setminus (-1,0] ^2}$, and the following Poisson problem:
\begin{align}\label{eq:poisson}
  -\Delta u &= 0,\quad \textrm{in } \Omega, \\
  u &= g_D, \quad \textrm{on } \partial\Omega, 
\end{align}
where $g_D$ corresponds to the Dirichlet trace of the analytical function in polar coordinates ${u(r,\theta)=r^\alpha \sin(\alpha \, \theta)}$, with ${\theta=3\pi/2}$ for our case. This particular problem is known as reentrant corner problem~\cite{ mitchell2013collection}, where the solution has a singularity at the corner, and its solution belongs to ${H^{1+\alpha-\epsilon}}$, ${\forall \epsilon>0}$~\cite{ oden1995parallel}.  The dG variational formulation we use in problem~\eqref{eq:poisson} is the formulation~\eqref{eq:gen_vf} with $K = 1$, ${\Vb = (0, 0)^T}$, ${\sigma=0}$ and ${f=0}$. Figure~\ref{fig:L-shape-43} shows the convergence plots for ${\alpha=2/3}$. A uniform refinement strategy cannot deal with the corner singularity achieving similar convergence rates for increasing polynomial orders.  However, the adaptive stabilized methodology overcomes this limitation, recovering optimal convergence rates on both $L^2$ and $V_h$ error norms.

\begin{figure}[ht!]
  \begin{center}
  	\hspace{0.08\textwidth}
  	\begin{subfigure}{0.325\textwidth}
  		\centering
  		\includegraphics[width=\textwidth]{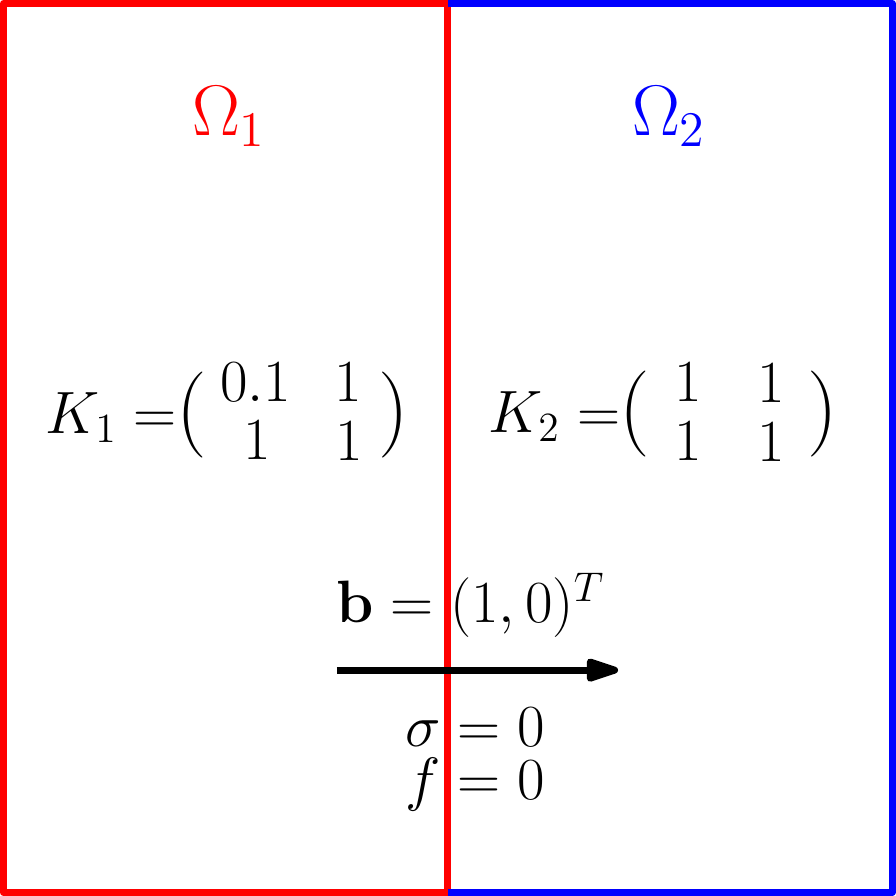} 
  		\label{figure5:sketch}
  		\caption{Problem sketch}
  		\vspace{0.05\textwidth}
  	\end{subfigure}
    \hspace{0.085\textwidth}
	\begin{subfigure}{0.48\textwidth}
  		\centering
  		\includegraphics[width=\textwidth]{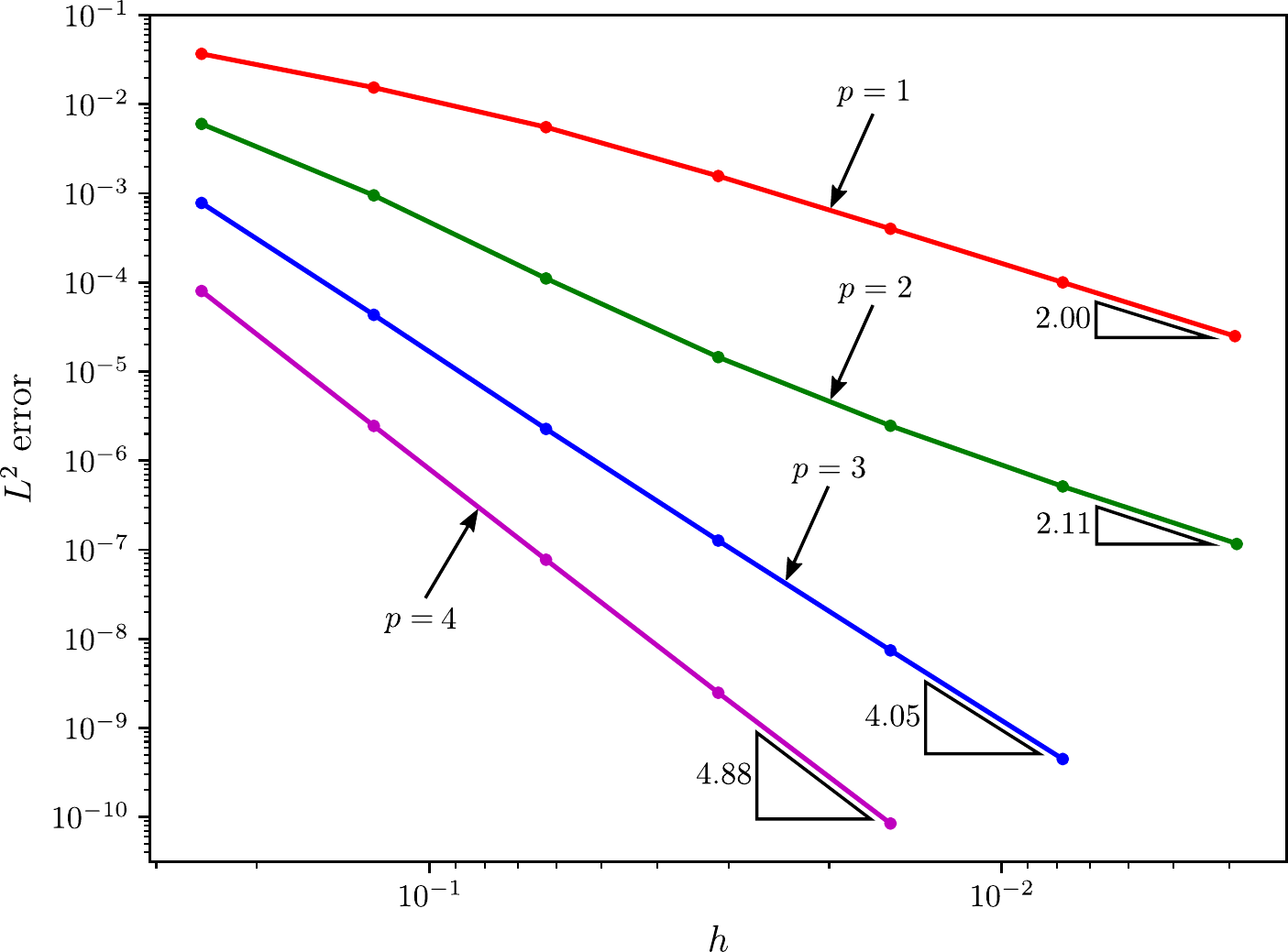} 
  		\label{figure5:L2}
  		\vspace{-0.05\textwidth}
  		\caption{Convergence plot in $L^2$ error norm}
  		\vspace{0.05\textwidth}
  	\end{subfigure}
    \begin{subfigure}{0.48\textwidth}
      \centering
      \includegraphics[width=\textwidth]{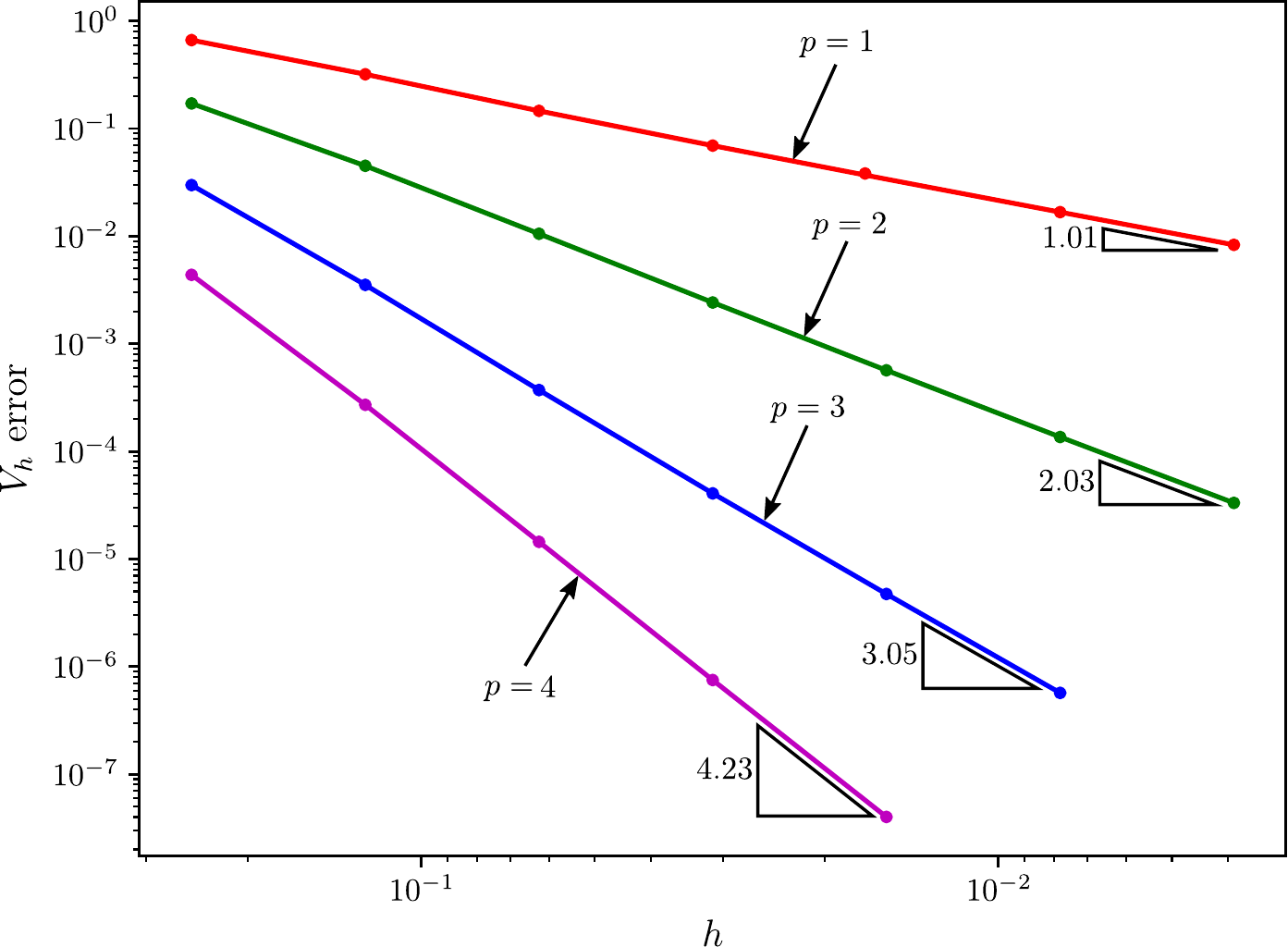} 
      \label{figure5:Vh}
      \vspace{-0.05\textwidth}
      \caption{Convergence plot in $V_h$ error norm}
      \vspace{0.05\textwidth}
    \end{subfigure}
	\hspace{0.02\textwidth}
    \begin{subfigure}{0.48\textwidth}
      \centering
      \includegraphics[width=\textwidth]{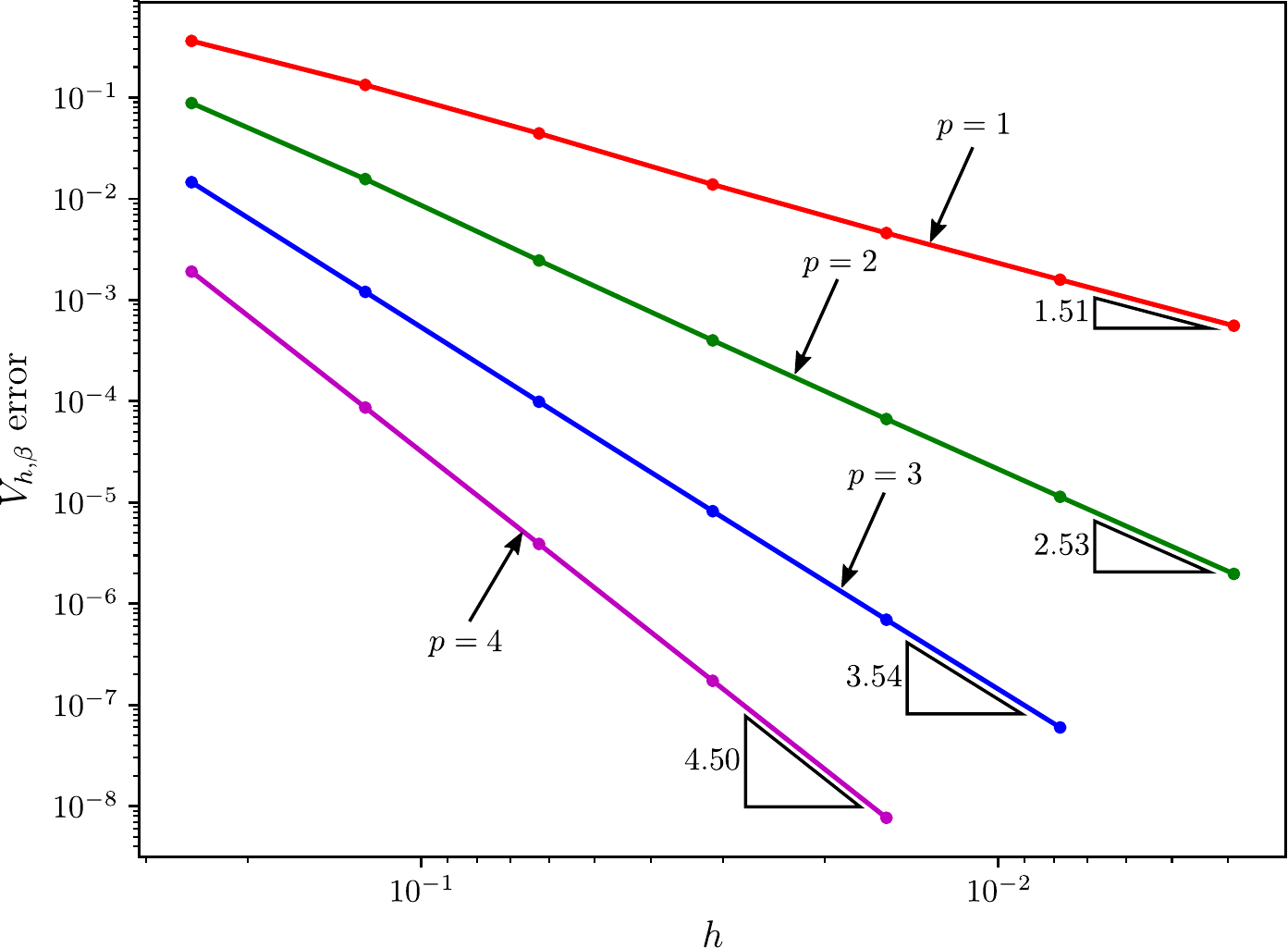} 
      \label{figure5:Vh_beta}
      \vspace{-0.05\textwidth}
      \caption{Convergence plot in $V_{h,\beta}$ error norm}
      \vspace{0.05\textwidth}
    \end{subfigure}
  \end{center}
  \caption{Heterogeneous diffusion problem.}
  \label{fig:hetero}
\end{figure}

\subsection{Heterogeneous diffusion}

We solve an advection-diffusion problem with discontinuous diffusion coefficients, based on a test from~\cite{ burman2006domain}; this problem shows the method performance with heterogeneous diffusion. We split the domain in two subdomains, $\Omega_1 = [0, \frac{1}{2}] \times [0,1]$ and $\Omega_2 = [\frac{1}{2}, 1] \times [0, 1]$ and use a constant diffusivity tensor in each subdomain 
\begin{align*}
  K_i(x,y)=\begin{pmatrix}
    \epsilon_i(x) & 0 \\
    0 & 1.0
  \end{pmatrix}
\end{align*}
where $\epsilon_i(x)$ represent discontinuous values across the interface $x=\frac{1}{2}$. We set $\epsilon_1(x)=1\times 10^{-2}$ and $\epsilon_2(x)=1.0$ with $\Vb = (1,0)^T$, ${\sigma=0}$ and ${f=0}$.  For this parameter choice, the exact solution is an exponential with respect to the $x$-coordinate (i.e., independent of the $y$-coordinate). At the interface, the solution satisfies the following conditions:
\begin{align*}
  \displaystyle \lim_{x\rightarrow\frac{1}{2}^-}u(x,y)&=\displaystyle \lim_{x\rightarrow\frac{1}{2}^+}u(x,y), \\
  \displaystyle \lim_{x\rightarrow\frac{1}{2}^-}-\epsilon(x) \partial_x u(x,y)&=\displaystyle \lim_{x\rightarrow\frac{1}{2}^+}-\epsilon(x) \partial_x u(x,y).
\end{align*} 

We set $u(0,y)=0$, $u(1,y)=1$, and by consequence of the matching conditions, we obtain
\begin{align*}
  u\left(\frac{1}{2},y\right)=\left[ \frac{u(0,y)\exp\left(\frac{1}{2\epsilon_1}\right)}{1-\exp\left(\frac{1}{2\epsilon_1}\right)} +\frac{u(1,y)}{1-\exp\left(\frac{1}{2\epsilon_2}\right)} \right]\left[ \frac{\exp\left(\frac{1}{2\epsilon_1}\right)}{1-\exp\left(\frac{1}{2\epsilon_1}\right)} +\frac{1}{1-\exp\left(\frac{1}{2\epsilon_2}\right)} \right]^{-1}.
\end{align*} 
Thus, the exact solution in each subdomain becomes:
\begin{align*}
  u_1(x,y)=&\frac{u\left(\frac{1}{2},y\right)-\exp\left(\frac{1}{2\epsilon_1}\right)u(0,y) + \left[u(0,y) - u\left(\frac{1}{2},y\right)\right]\exp\left(\frac{x}{\epsilon_1}\right)}{1-\exp\left(\frac{1}{2\epsilon_1}\right)}\\
  u_1(x,y)=&\frac{u\left(1,y\right)-\exp\left(\frac{1}{2\epsilon_2}\right)u(\frac{1}{2},y) + \left[u(\frac{1}{2},y) - u\left(1,y\right)\right]\exp\left(\frac{x-\frac{1}{2}}{\epsilon_2}\right)}{1-\exp\left(\frac{1}{2\epsilon_2}\right)}
\end{align*}
Figure~\ref{fig:hetero} shows the convergence rates that using both trial and test space functions of the same polynomial degree deliver, for $p=1,2,3$. As expected, we recover the same convergence rates as the original dG scheme: $h^p$ in the $V_h$ error norm and $h^{p+1/2}$ in the $V_{h, \beta}$ error norm.

\begin{figure}[ht!]
	\centering
	\includegraphics[width=0.35\textwidth]{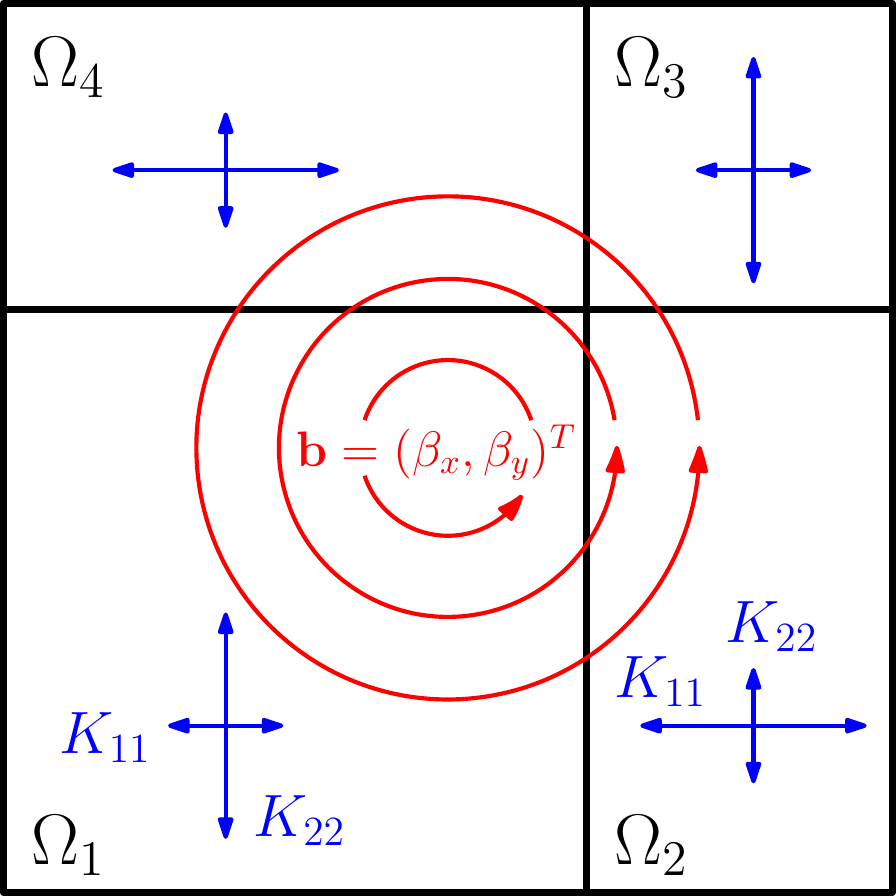} 
	\caption{Anisotropic diffusion problem sketch. Counterclockwise advection field.}
	\label{fig:sketch}
\end{figure}

\subsection{Anisotropic with high-contrast diffusion}

In this example, we consider an anisotropic problem, following~\cite{ ern2009}. We consider the unit square $\Omega = [0,1]\times[0,1]$ split into four subdomains: $\Omega_1 = \left[0,\frac{2}{3}\right]\times\left[0,\frac{2}{3}\right]$, $\Omega_2 = \left[\frac{2}{3},1\right]\times\left[0,\frac{2}{3}\right]$, $\Omega_3 = \left[\frac{2}{3},1\right]\times\left[\frac{2}{3},1\right]$ and $\Omega_4 = \left[0,\frac{2}{3}\right]\times\left[\frac{2}{3},1\right]$. The diffusivity tensor takes different values in each subdomain , see Figure~\ref{fig:sketch}:
\begin{align*}
K_i(x,y)&=\begin{pmatrix}
10^{-6} & 0 \\
0 & 1.0
\end{pmatrix} &&\text{for } i=1,3\ \ \forall (x,y) \in \Omega_1, \Omega_3, \\
K_i(x,y)&=\begin{pmatrix}
1.0 & 0 \\
0 & 10^{-6}
\end{pmatrix} &&\text{for } i=2,4\ \  \forall (x,y) \in \Omega_2, \Omega_4.
\end{align*}

\begin{figure}[ht!]
  \begin{center}
    \begin{subfigure}{0.45\textwidth}
      \centering
      \includegraphics[width=\textwidth]{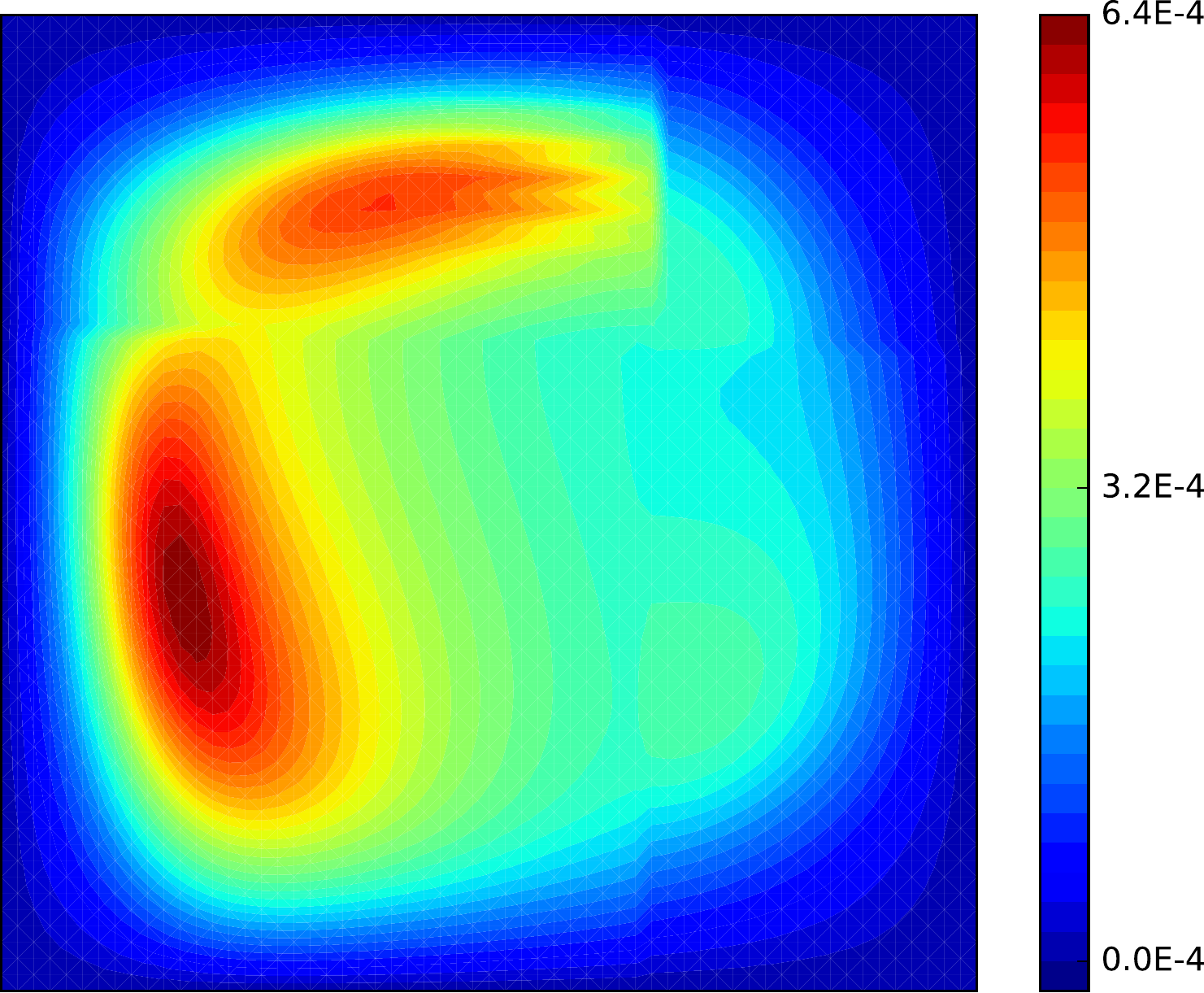} 
      \label{anisotropic-stephansen:counter}
      \caption{Counterclockwise advection field}
    \end{subfigure}
	\hspace{0.05\textwidth}
    \begin{subfigure}{0.45\textwidth}
      \centering
      \includegraphics[width=\textwidth]{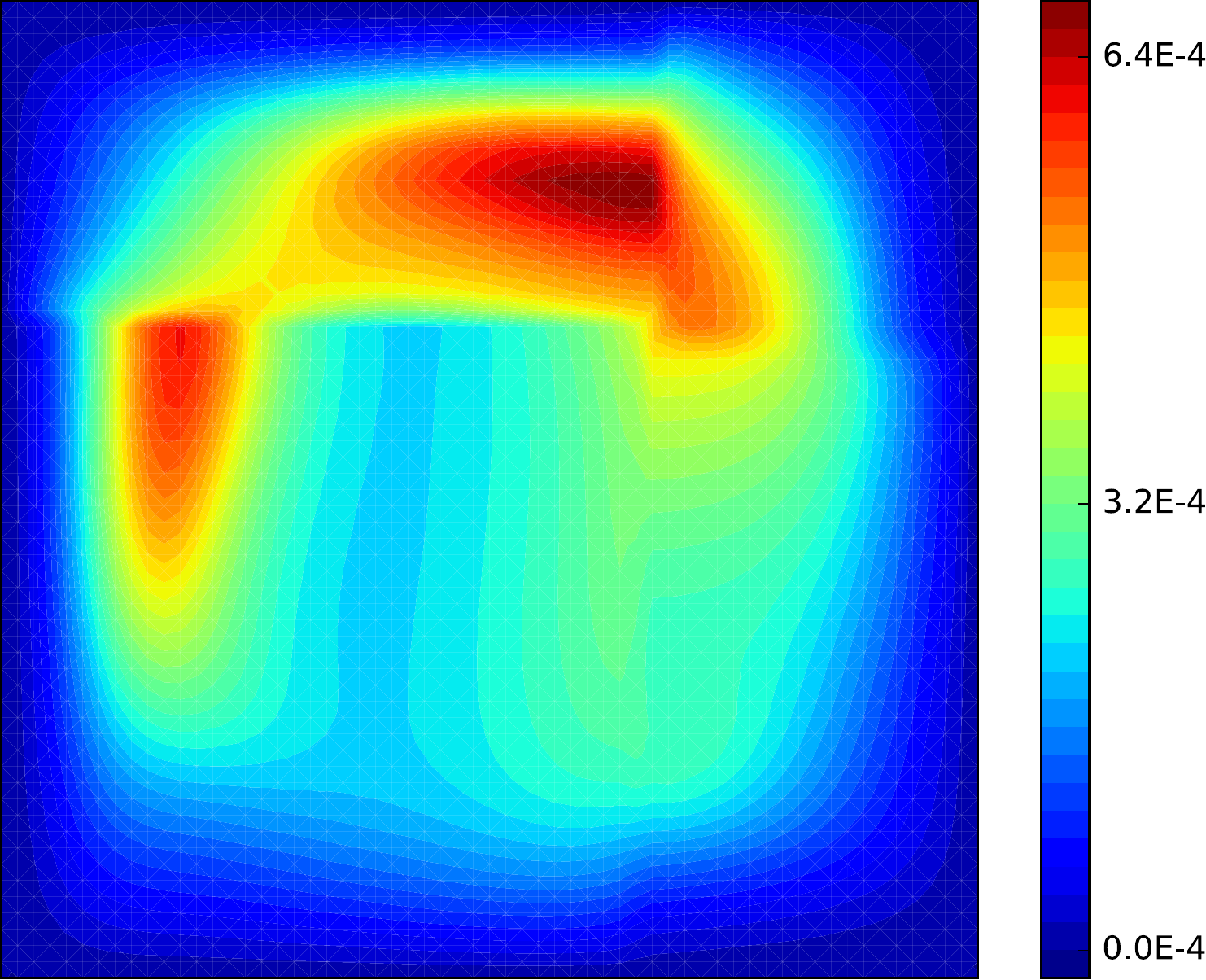} 
      \label{anisotropic-stephansen:clock}
      \caption{Clockwise advection field}
    \end{subfigure}
  \end{center}
  \caption{Anisotropic diffusion problem. Uniform mesh (25.3~K~DOFs).}
  \label{fig:anisotropic-stephansen}
\end{figure}
      
\begin{figure}[t!]
	\begin{center}
		\begin{subfigure}{0.45\textwidth}
			\centering
			\includegraphics[width=\textwidth]{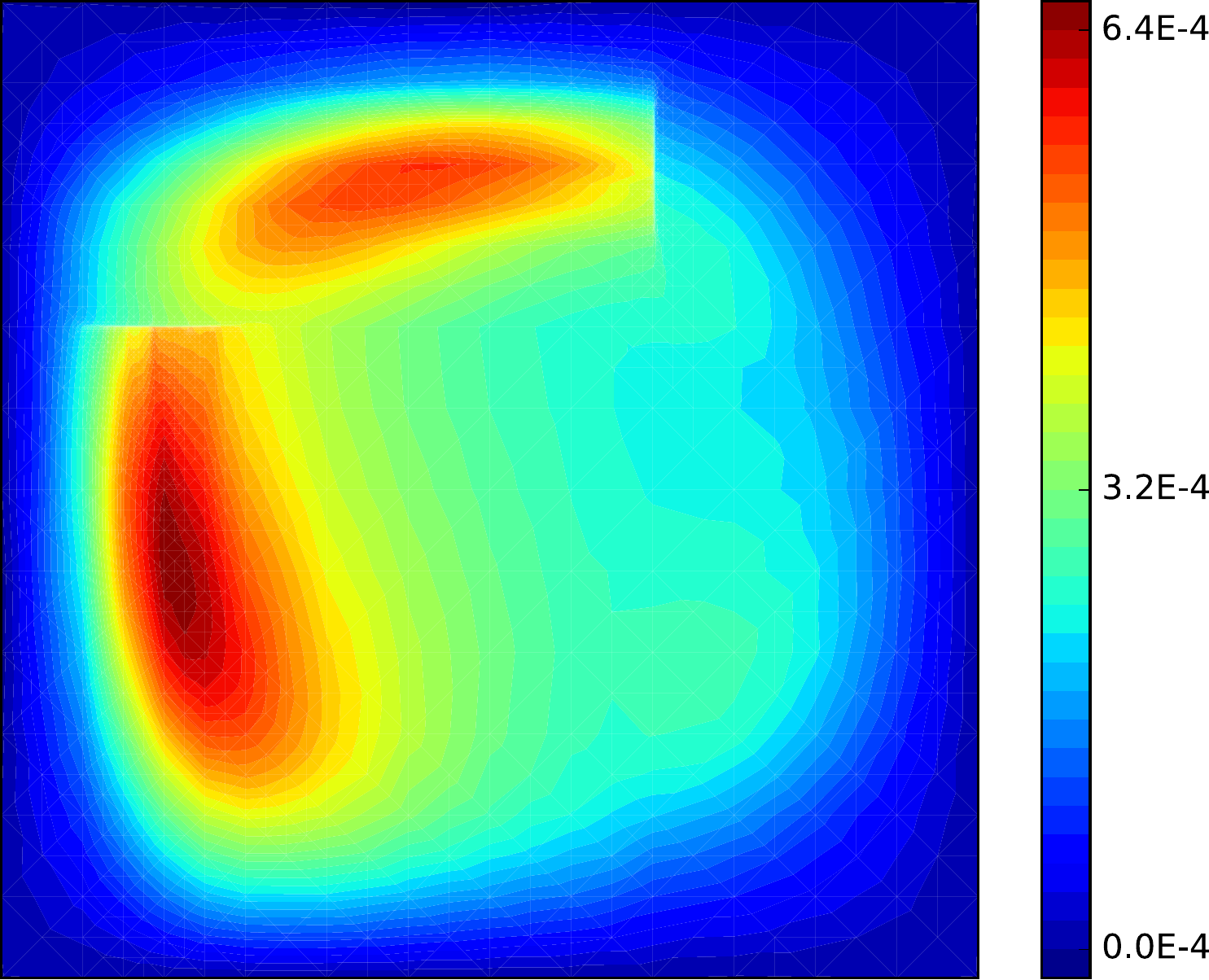} 
			\label{anisotropic-adaptive:counter}
			\caption{Computed solution and refined mesh, counterclockwise advection field (Level 12: 23.7~K~DOFs)}
		\end{subfigure}
		\hspace{0.05\textwidth}
		\begin{subfigure}{0.45\textwidth}
			\centering
			\includegraphics[width=\textwidth]{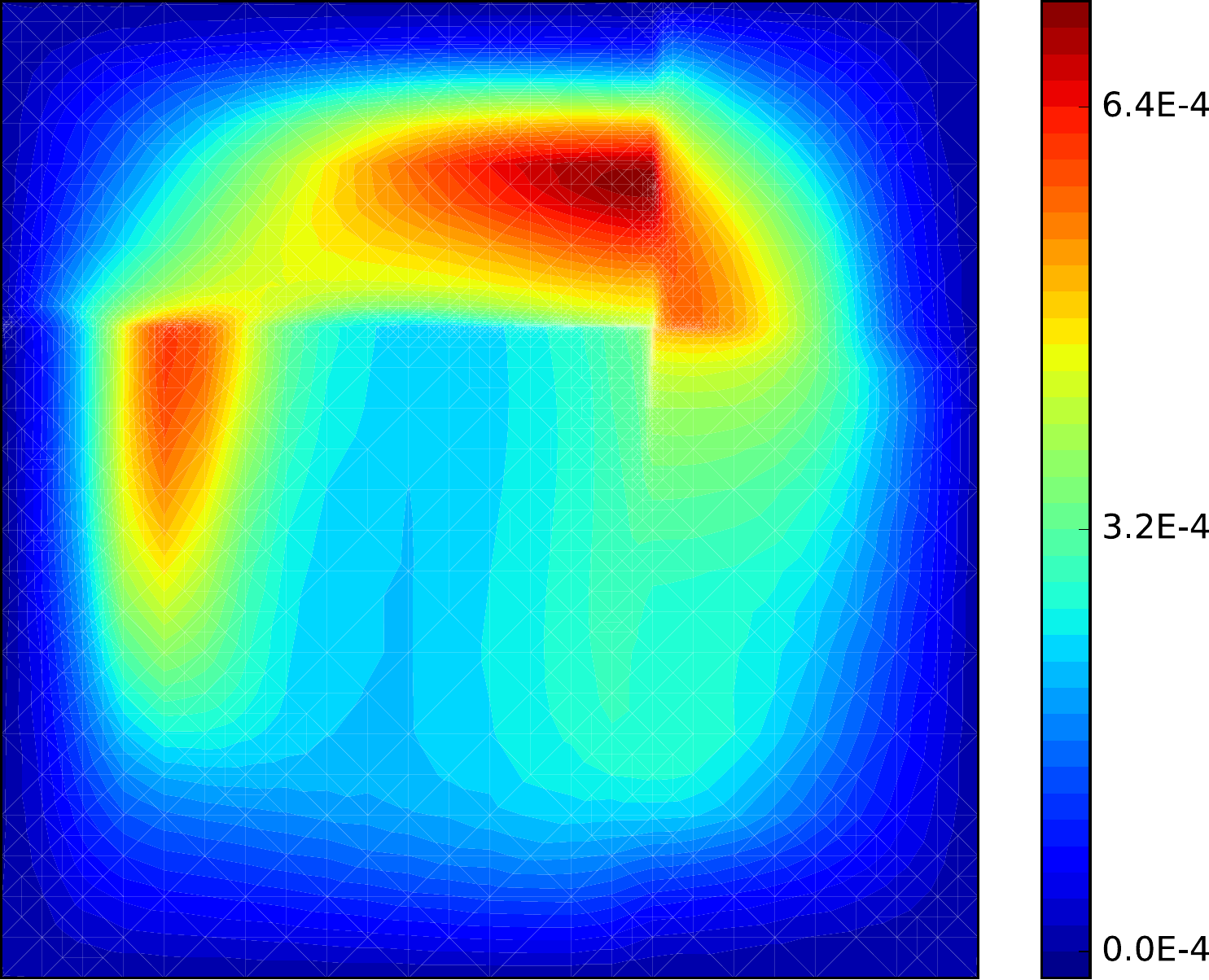} 
			\label{anisotropic-adaptive:clock}
			\caption{Computed solution and refined mesh, clockwise advection field (Level 10: 22.3~K~DOFs)}
		\end{subfigure}
	\end{center}
	\caption{Anisotropic diffusion problem. Adaptive refinement.}
	\label{fig:anisotropic-adaptive}
\end{figure}

The advection field is solenoidal ${\Vb = (\beta_x, \beta_y)^T}$, with ${\beta_x=40x(2y-1)(x-1)}$ and ${\beta_y=-40y(2x-1)(y-1)}$ for the counterclockwise case, and ${\Vb = -(\beta_x, \beta_y)^T}$ for the clockwise case. Unlike the previous example, the advective field is neither constant nor orthogonal to the discontinuities in $K$. However, its orientation is still along the direction of increasing diffusivity; for that reason, internal layers develop in the solution. The forcing term is ${f(x,y)=10^{-2}\exp(-(r-0.35)^2/0.005)}$ with ${r^2 = (x-0.5)^2+(y-0.5)^2}$, corresponding to a Gaussian hill with center at $r=0.35$. Finally, we set $\sigma=1$ for the reaction term, and $g_D=0$ on $\partial\Omega$ for the boundary condition. We consider two subcases, the first on a quasi-uniform mesh with ${h=0.024}$, conforming to the discontinuities of $K$, and the second through an adaptive scheme starting from a uniform triangular mesh (${h=0.177}$). We solve both cases using the same polynomial degree ${p=1}$ for trial and test spaces.  Figure~\ref{fig:anisotropic-stephansen} shows results obtained for a uniform mesh, while Figure~\ref{fig:anisotropic-adaptive} shows results for adaptive mesh refinement with fewer degrees of freedom than in the uniform mesh case. We use the SWIP formulation such that the jump penalty considers the principal directions of the diffusivity. Thus, SWIP avoids overshoots and undershoots near the material interfaces, which is not the case for standard interior penalty schemes~\cite{ ern2009}. Finally, the adaptive strategy concentrates refinement at the inner layer but without losing the approximation quality.

\begin{figure}[t!]
  \begin{center}
    \begin{subfigure}{0.38\textwidth}
      \centering
      \includegraphics[width=\textwidth]{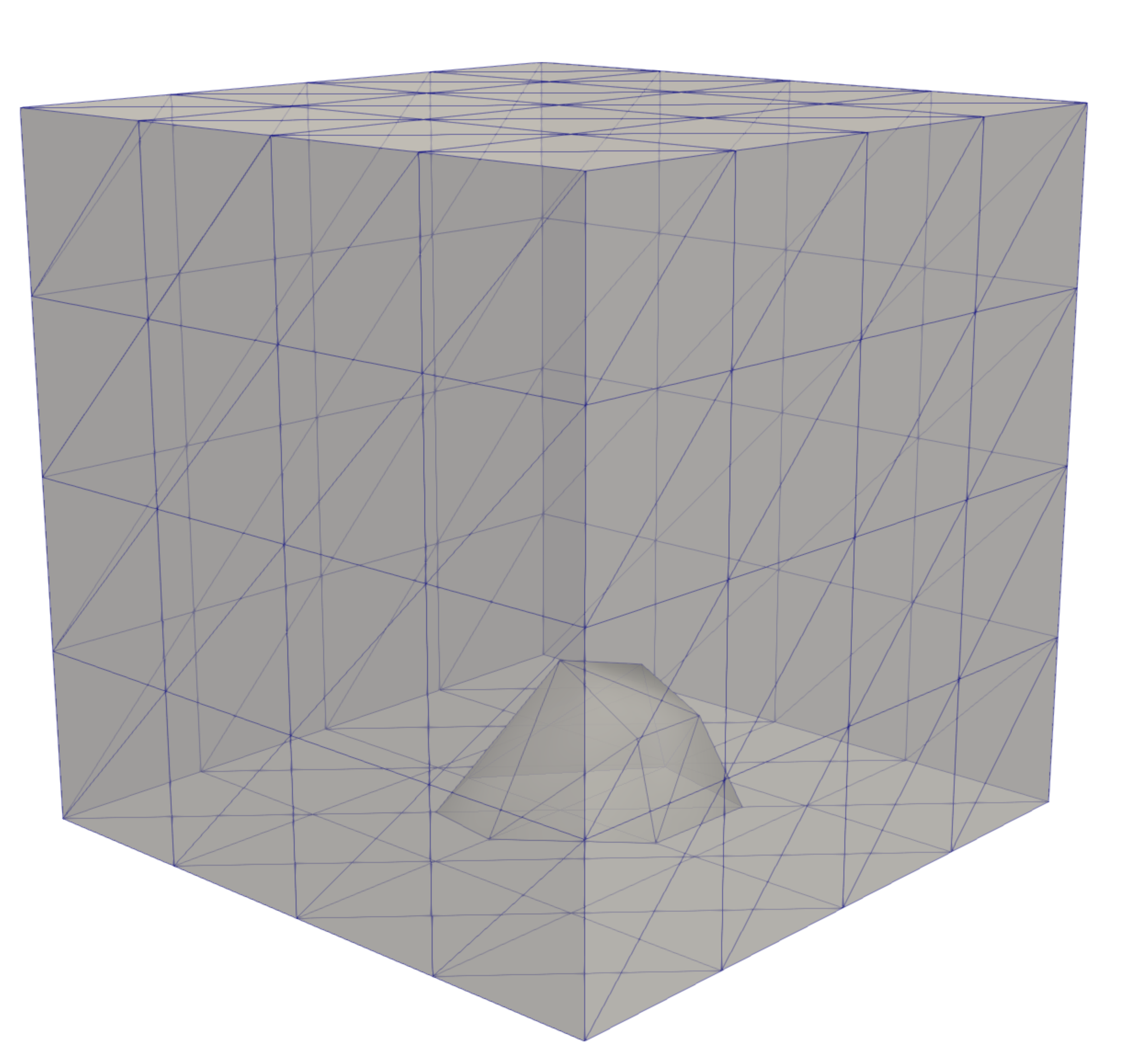} 
      \caption{Level 0. Inner contour for $u=0.5$.}
      \label{3D:mesh}
    \end{subfigure}
    \hspace{0.1\textwidth}
    \begin{subfigure}{0.38\textwidth}
      \centering
      \includegraphics[width=\textwidth]{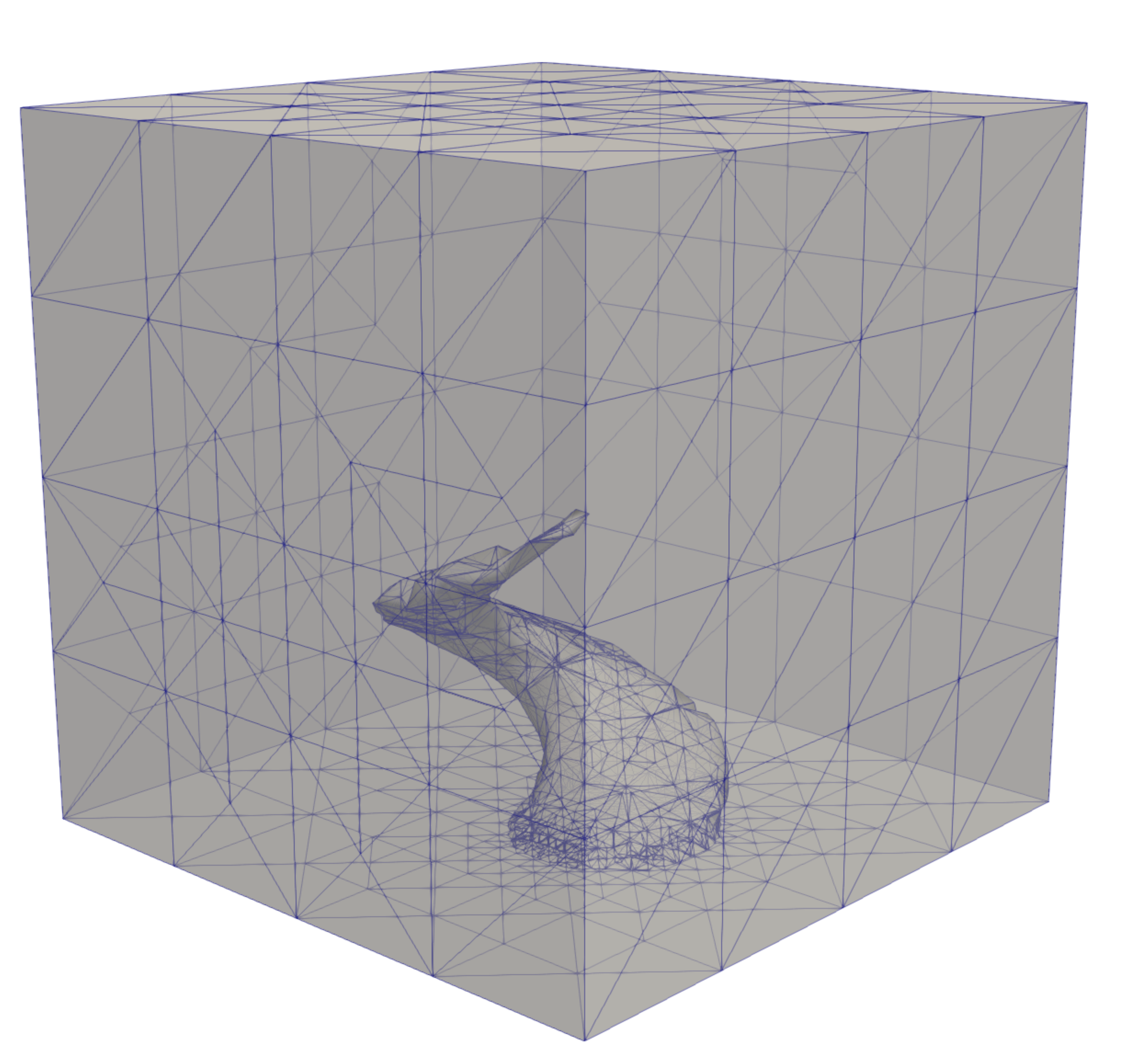} 
      \caption{Level 5. Inner contour for $u=1$.}
      \label{3D:sol}
    \end{subfigure}
    \begin{subfigure}{0.38\textwidth}
      \centering
      \includegraphics[width=\textwidth]{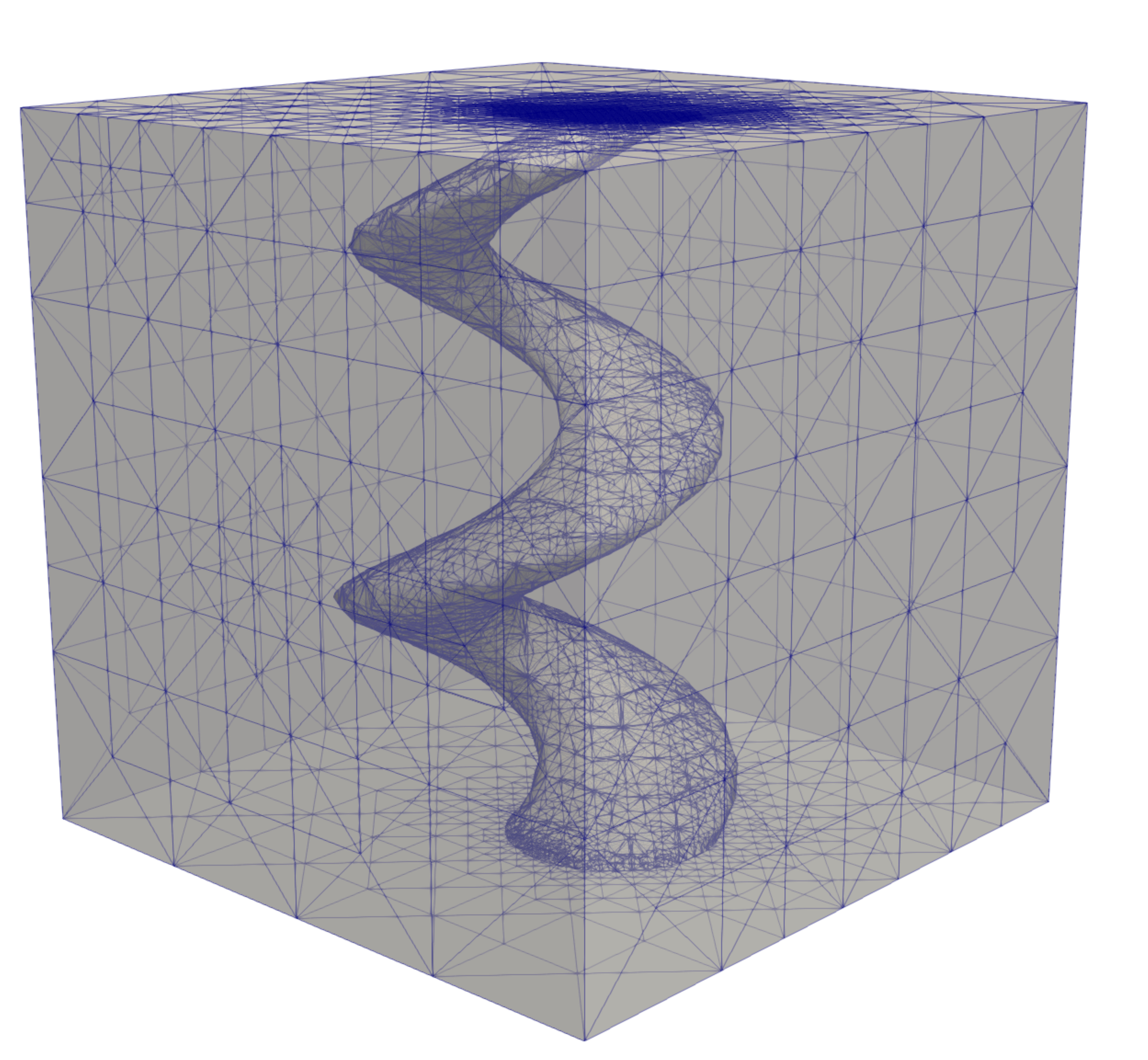} 
      \caption{Level 9. Inner contour for $u=1$.}
      \label{3D:bottom}
    \end{subfigure}
    \hspace{0.1\textwidth}
    \begin{subfigure}{0.38\textwidth}
      \centering
      \includegraphics[width=\textwidth]{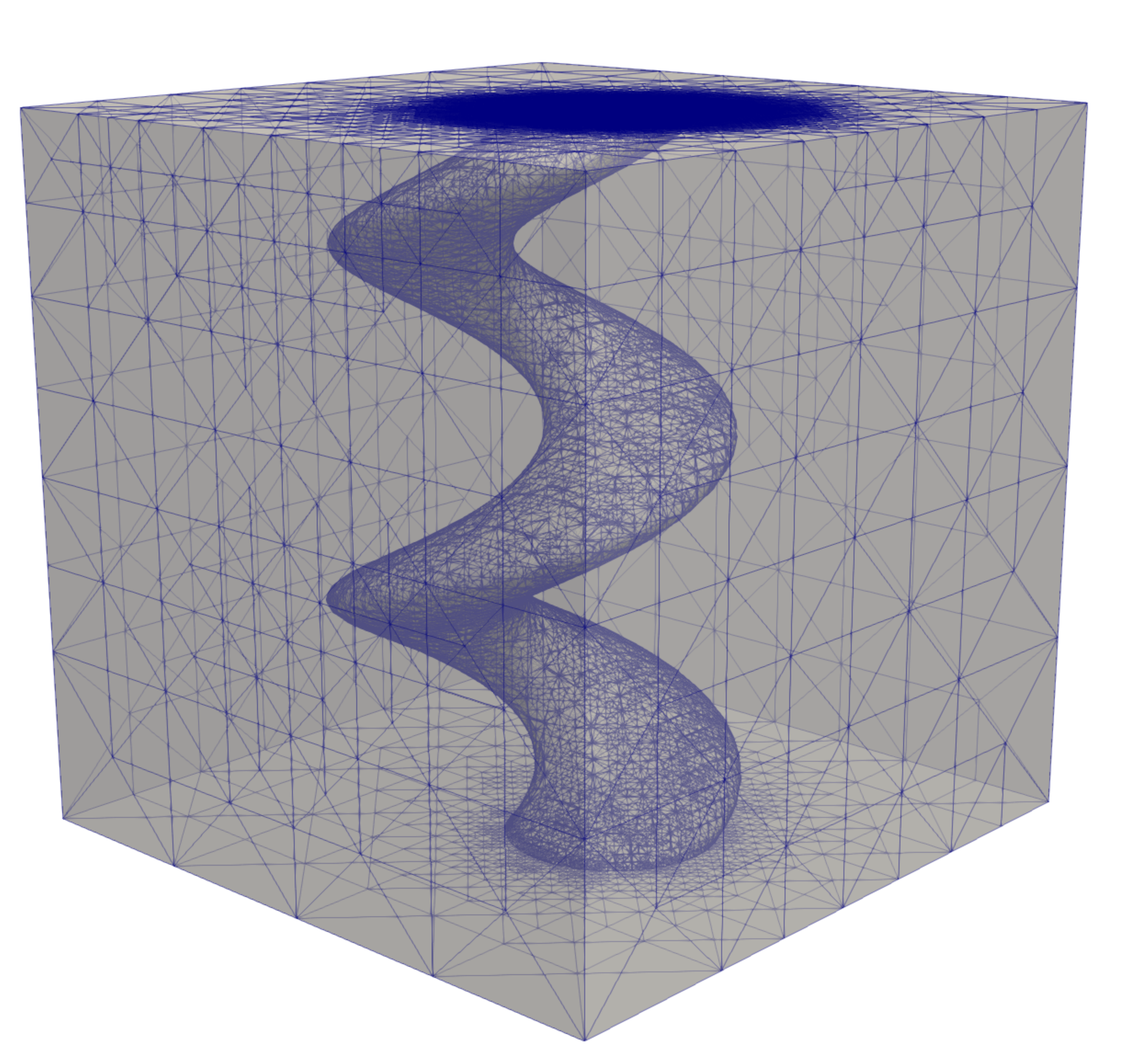} 
      \caption{Level 13. Inner contour for $u=1$.}
      \label{3D:top}
    \end{subfigure}
  \end{center}
  \caption{3D advection-dominated diffusion problem. Adaptive mesh evolution. }
  \label{fig:3DMesh_evolution}
\end{figure}

\subsection{3D advection-dominated diffusion}

\begin{figure}[t!]
  \centering
  \includegraphics[width=0.5\textwidth]{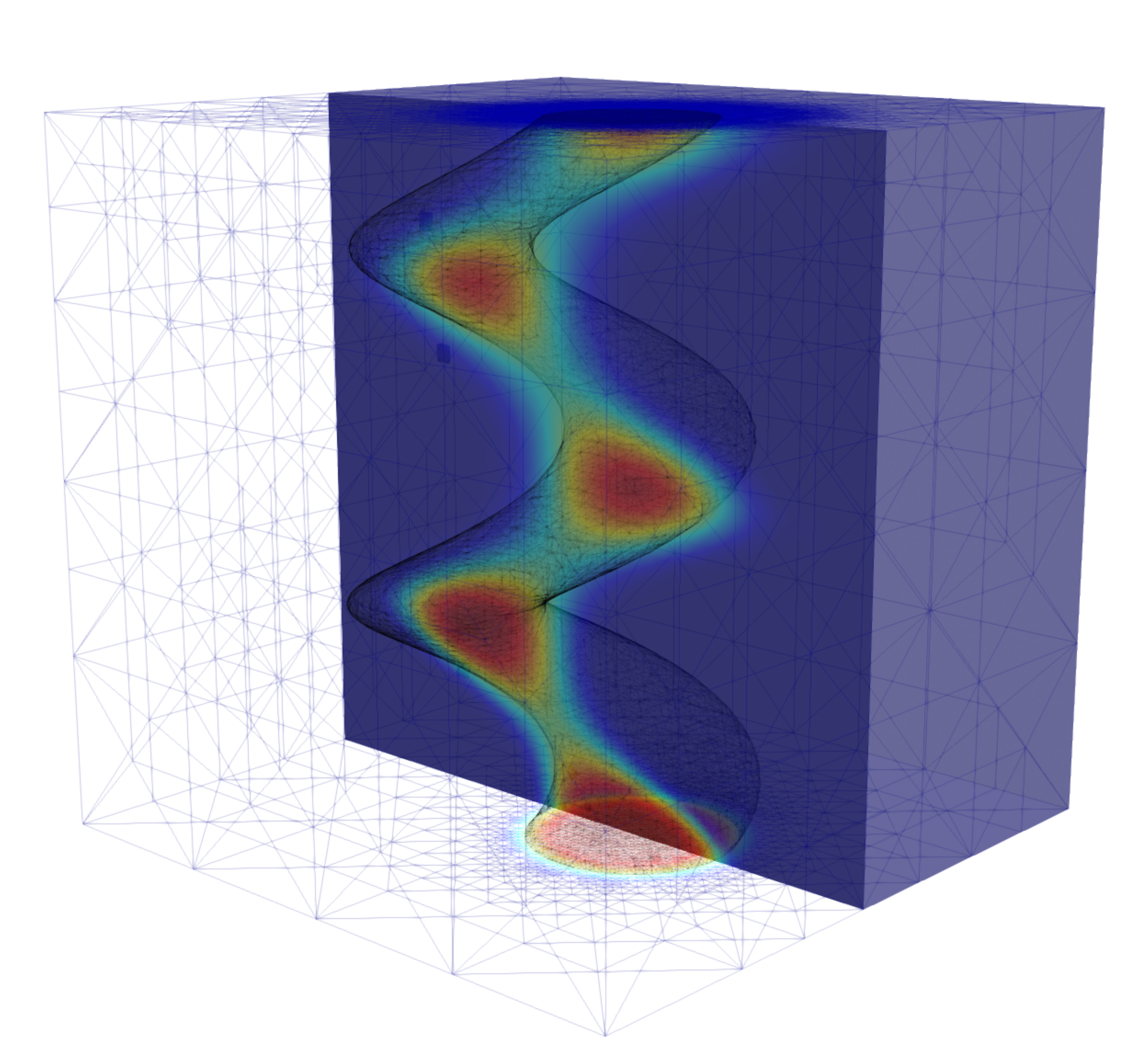} 
  \caption{3D advection-dominated diffusion problem. Slice over $y=0.5$. Level 13: 4'858,125~DOFs. }
  \label{fig:3D}
\end{figure}

As a last example, we consider a 3D advection-diffusion problem in the unit cube ${\Omega = (0 \, , \, 1)^3 \subset \R^3}$. We set the source term ${f=0}$, the diffusion ${K=10^{-3}}$, the spiral-type advective field ${\Vb=(\beta_x,\beta_y,\beta_z)^T =  (-0.15\,\sin(4\pi z) \, , \, 0.15\,\cos(4\pi z) \, , \, 1)^T}$,  and the inflow boundary datum $g$ as
\begin{equation*}
  g =
  \left\{
    \renewcommand{\arraystretch}{1.4}
    \begin{array}{ll}
      1+\tanh\left[M \left( 0.15^2 - \left(x- 0.6 \right)^2 - \left(y- 0.5 \right)^2 \right) \right] & \textrm{on } z=0, \\
      0 & \textrm{elsewhere on} \ \Gamma, \\
    \end{array}
  \right.
\end{equation*} 
These parameters produce a solution that presents an inner layer in the bottom of the unit cube (${z=0}$) for values of $M\gg 1$ in the inflow boundary datum. On the other hand, at the top of the unit cube (${z=1}$), the solution exhibits a boundary layer due to the advection-dominant regime and the boundary condition imposed on that surface. Figure~\ref{fig:3DMesh_evolution} shows the evolution of the 3D mesh as the refinement strategy progresses. Figure~\ref{3D:mesh} displays the initial progress of the spiral shape of the inner layer inside the unit cube. Figures~\ref{3D:sol},~\ref{3D:bottom}, and~\ref{3D:top} reflect the refinement process to capture the sharp internal layer that the advective field induces as well as the boundary layer that appears at the outlet due to the small diffusion. Figure~\ref{fig:3D}, shoes the discrete solution for a refined mesh displaying the interior mesh as well as a cut of the solution.

\section{Conclusions}\label{ref:concl} 

In this paper, we describe an adaptive stabilized conforming finite element method that minimizes the residual on dual discontinuous Galerkin (dG) norms for advection-diffusion-reaction problems.  We demonstrate the performance and robustness of the method by analyzing in detail several challenging diffusivity distributions. In particular, we tackle problems with highly heterogeneous and anisotropic coefficient distributions as well as advection-dominated cases. The method computes smooth solutions that minimize a residual on a dual dG norm. This residual minimization problem leads to an equivalent saddle-point problem, which we solve explicitly for a discretely continuous approximation and an efficient and reliable error estimate that can guide adaptive mesh refinement. In the context of advection-diffusion-reaction problems, this method recovers the optimal convergence rates for $h$-adaptive schemes in the dG norm. Besides, the method captures sharp boundary and internal layers as well as overcomes the classical overshooting and undershooting problems. In summary, we use several challenging numerical examples to show the robustness of the method under a wide range of scenarios. 

\section{Acknowledgments}

The CSIRO Deep Earth Imaging Future Science Platform Postgraduate Top-Up Scholarship supports the work of RJC and this is gratefully acknowledged. This publication was made possible in part by the CSIRO Professorial Chair in Computational Geoscience at Curtin University and the Deep Earth Imaging Enterprise Future Science Platforms of the Commonwealth Scientific Industrial Research Organisation, CSIRO, of Australia. This project has also received funding from the European Union's Horizon 2020 research and innovation programme under the Marie Sklodowska-Curie grant agreement No 777778 (MATHROCKS). At Curtin University, The Curtin Corrosion Centre, the Curtin Institute for Computation, and The Institute for Geoscience Research (TIGeR) kindly provide continuing support. 

\bibliography{mybibfile}

\begin{thebibliography}{57}
\expandafter\ifx\csname natexlab\endcsname\relax\def\natexlab#1{#1}\fi
\providecommand{\url}[1]{\texttt{#1}}
\providecommand{\href}[2]{#2}
\providecommand{\path}[1]{#1}
\providecommand{\DOIprefix}{doi:}
\providecommand{\ArXivprefix}{arXiv:}
\providecommand{\URLprefix}{URL: }
\providecommand{\Pubmedprefix}{pmid:}
\providecommand{\doi}[1]{\href{http://dx.doi.org/#1}{\path{#1}}}
\providecommand{\Pubmed}[1]{\href{pmid:#1}{\path{#1}}}
\providecommand{\bibinfo}[2]{#2}
\ifx\xfnm\relax \def\xfnm[#1]{\unskip,\space#1}\fi
\bibitem[{Calo et~al.(2020)Calo, Ern, Muga, and Rojas}]{calo2019}
\bibinfo{author}{V.~M. Calo}, \bibinfo{author}{A.~Ern},
  \bibinfo{author}{I.~Muga}, \bibinfo{author}{S.~Rojas},
\newblock \bibinfo{title}{{An adaptive stabilized conforming finite element
  method via residual minimization on dual discontinuous Galerkin norms}},
\newblock \bibinfo{journal}{Computer Methods in Applied Mechanics and
  Engineering} \bibinfo{volume}{363} (\bibinfo{year}{2020})
  \bibinfo{pages}{112891}.
\bibitem[{Ewing and Wang(2001)}]{ewing2001}
\bibinfo{author}{R.~E. Ewing}, \bibinfo{author}{H.~Wang},
\newblock \bibinfo{title}{A summary of numerical methods for time-dependent
  advection-dominated partial differential equations},
\newblock \bibinfo{journal}{Journal of Computational and Applied Mathematics}
  \bibinfo{volume}{128} (\bibinfo{year}{2001}) \bibinfo{pages}{423--445}.
\bibitem[{Calo et~al.(2014)Calo, Efendiev, Galvis, and Ghommem}]{Calo2014}
\bibinfo{author}{V.~M. Calo}, \bibinfo{author}{Y.~Efendiev},
  \bibinfo{author}{J.~Galvis}, \bibinfo{author}{M.~Ghommem},
\newblock \bibinfo{title}{{Multiscale empirical interpolation for solving
  nonlinear PDEs}},
\newblock \bibinfo{journal}{Journal of Computational Physics}
  \bibinfo{volume}{278} (\bibinfo{year}{2014}) \bibinfo{pages}{204 -- 220}.
\bibitem[{Ern et~al.(2009)Ern, Stephansen, and Zunino}]{ern2009}
\bibinfo{author}{A.~Ern}, \bibinfo{author}{A.~F. Stephansen},
  \bibinfo{author}{P.~Zunino},
\newblock \bibinfo{title}{{A discontinuous Galerkin method with weighted
  averages for advection--diffusion equations with locally small and
  anisotropic diffusivity}},
\newblock \bibinfo{journal}{IMA Journal of Numerical Analysis}
  \bibinfo{volume}{29} (\bibinfo{year}{2009}) \bibinfo{pages}{235--256}.
\bibitem[{Hossain et~al.(2012)Hossain, Hossainy, Bazilevs, Calo, and
  Hughes}]{Hossain2012}
\bibinfo{author}{S.~S. Hossain}, \bibinfo{author}{S.~F.~A. Hossainy},
  \bibinfo{author}{Y.~Bazilevs}, \bibinfo{author}{V.~M. Calo},
  \bibinfo{author}{T.~J.~R. Hughes},
\newblock \bibinfo{title}{Mathematical modeling of coupled drug and
  drug-encapsulated nanoparticle transport in patient-specific coronary artery
  walls},
\newblock \bibinfo{journal}{Computational Mechanics} \bibinfo{volume}{49}
  (\bibinfo{year}{2012}) \bibinfo{pages}{213--242}.
\bibitem[{Calo et~al.(2008)Calo, Brasher, Bazilevs, and Hughes}]{Calo2008}
\bibinfo{author}{V.~Calo}, \bibinfo{author}{N.~Brasher},
  \bibinfo{author}{Y.~Bazilevs}, \bibinfo{author}{T.~Hughes},
\newblock \bibinfo{title}{Multiphysics model for blood flow and drug transport
  with application to patient-specific coronary artery flow},
\newblock \bibinfo{journal}{Computational Mechanics} \bibinfo{volume}{43}
  (\bibinfo{year}{2008}) \bibinfo{pages}{161--177}.
\bibitem[{Bazilevs et~al.(2007)Bazilevs, Calo, Tezduyar, and
  Hughes}]{Bazilevs2007}
\bibinfo{author}{Y.~Bazilevs}, \bibinfo{author}{V.~M. Calo},
  \bibinfo{author}{T.~E. Tezduyar}, \bibinfo{author}{T.~J.~R. Hughes},
\newblock \bibinfo{title}{{\textit{YZ}}$\beta$ discontinuity capturing for
  advection-dominated processes with application to arterial drug delivery},
\newblock \bibinfo{journal}{International Journal for Numerical Methods in
  Fluids} \bibinfo{volume}{54} (\bibinfo{year}{2007})
  \bibinfo{pages}{593--608}.
\bibitem[{Calo et~al.(2016)Calo, Efendiev, Galvis, and Li}]{Calo2016}
\bibinfo{author}{V.~M. Calo}, \bibinfo{author}{Y.~Efendiev},
  \bibinfo{author}{J.~Galvis}, \bibinfo{author}{G.~Li},
\newblock \bibinfo{title}{Randomized oversampling for generalized multiscale
  finite element methods},
\newblock \bibinfo{journal}{Multiscale Modeling \& Simulation}
  \bibinfo{volume}{14} (\bibinfo{year}{2016}) \bibinfo{pages}{482--501}.
\bibitem[{Galvis et~al.(2018)Galvis, Chung, Efendiev, and Leung}]{galvis2018}
\bibinfo{author}{J.~Galvis}, \bibinfo{author}{E.~Chung},
  \bibinfo{author}{Y.~Efendiev}, \bibinfo{author}{W.~Leung},
\newblock \bibinfo{title}{On overlapping domain decomposition methods for
  high-contrast multiscale problems},
\newblock in: \bibinfo{editor}{B.~P. et~al.} (Ed.), \bibinfo{booktitle}{Domain
  Decomposition Methods in Science and Engineering XXIV. DD 2017}, volume
  \bibinfo{volume}{125} of \textit{\bibinfo{series}{Lecture Notes in
  Computational Science and Engineering}}, \bibinfo{publisher}{Springer},
  \bibinfo{year}{2018}.
\bibitem[{Calo et~al.(2011)Calo, Efendiev, and Galvis}]{Calo2011}
\bibinfo{author}{V.~Calo}, \bibinfo{author}{Y.~Efendiev},
  \bibinfo{author}{J.~Galvis},
\newblock \bibinfo{title}{A note on variational multiscale methods for
  high-contrast heterogeneous porous media flows with rough source terms},
\newblock \bibinfo{journal}{Advances in Water Resources} \bibinfo{volume}{34}
  (\bibinfo{year}{2011}) \bibinfo{pages}{1177 -- 1185}. \bibinfo{note}{New
  Computational Methods and Software Tools}.
\bibitem[{Brooks and Hughes(1982)}]{brooks1982}
\bibinfo{author}{A.~N. Brooks}, \bibinfo{author}{T.~J. Hughes},
\newblock \bibinfo{title}{{Streamline upwind/Petrov-Galerkin formulations for
  convection dominated flows with particular emphasis on the incompressible
  Navier-Stokes equations}},
\newblock \bibinfo{journal}{{Computer Methods in Applied Mechanics and
  Engineering}} \bibinfo{volume}{32} (\bibinfo{year}{1982})
  \bibinfo{pages}{199--259}.
\bibitem[{Hughes et~al.(1989)Hughes, Franca, and Hulbert}]{hughes1989new}
\bibinfo{author}{T.~J. Hughes}, \bibinfo{author}{L.~P. Franca},
  \bibinfo{author}{G.~M. Hulbert},
\newblock \bibinfo{title}{{A new finite element formulation for computational
  fluid dynamics: VIII. The Galerkin/least-squares method for
  advective-diffusive equations}},
\newblock \bibinfo{journal}{{Computer Methods in Applied Mechanics and
  Engineering}} \bibinfo{volume}{73} (\bibinfo{year}{1989})
  \bibinfo{pages}{173--189}.
\bibitem[{Hughes and Sangalli(2007)}]{hughes2007}
\bibinfo{author}{T.~J. Hughes}, \bibinfo{author}{G.~Sangalli},
\newblock \bibinfo{title}{{Variational multiscale analysis: the fine-scale
  Green’s function, projection, optimization, localization, and stabilized
  methods}},
\newblock \bibinfo{journal}{{SIAM Journal on Numerical Analysis}}
  \bibinfo{volume}{45} (\bibinfo{year}{2007}) \bibinfo{pages}{539--557}.
\bibitem[{Brezzi et~al.(1998)Brezzi, Marini, and Russo}]{brezzi1998}
\bibinfo{author}{F.~Brezzi}, \bibinfo{author}{D.~Marini},
  \bibinfo{author}{A.~Russo},
\newblock \bibinfo{title}{Applications of the pseudo residual-free bubbles to
  the stabilization of convection-diffusion problems},
\newblock \bibinfo{journal}{Computer Methods in Applied Mechanics and
  Engineering} \bibinfo{volume}{166} (\bibinfo{year}{1998})
  \bibinfo{pages}{51--63}.
\bibitem[{Reed and Hill(1973)}]{reed1973}
\bibinfo{author}{W.~H. Reed}, \bibinfo{author}{T.~Hill},
  \bibinfo{title}{{Triangular mesh methods for the neutron transport
  equation}}, \bibinfo{type}{Technical Report}, Los Alamos Scientific Lab., N.
  Mex.(USA), \bibinfo{year}{1973}.
\bibitem[{Lesaint and Raviart(1974)}]{lesaint1974}
\bibinfo{author}{P.~Lesaint}, \bibinfo{author}{P.-A. Raviart},
\newblock \bibinfo{title}{On a finite element method for solving the neutron
  transport equation},
\newblock \bibinfo{journal}{{Publications math{\'e}matiques et informatique de
  Renne}s}  (\bibinfo{year}{1974}) \bibinfo{pages}{1--40}.
\bibitem[{Johnson and Pitk{\"a}ranta(1986)}]{johnson1986}
\bibinfo{author}{C.~Johnson}, \bibinfo{author}{J.~Pitk{\"a}ranta},
\newblock \bibinfo{title}{{An analysis of the discontinuous Galerkin method for
  a scalar hyperbolic equation}},
\newblock \bibinfo{journal}{{Mathematics of Computation}} \bibinfo{volume}{46}
  (\bibinfo{year}{1986}) \bibinfo{pages}{1--26}.
\bibitem[{Cockburn et~al.(2012)Cockburn, Karniadakis, and Shu}]{cockburn2012}
\bibinfo{author}{B.~Cockburn}, \bibinfo{author}{G.~E. Karniadakis},
  \bibinfo{author}{C.-W. Shu}, \bibinfo{title}{{Discontinuous Galerkin methods:
  theory, computation and applications}}, volume~\bibinfo{volume}{11},
  \bibinfo{publisher}{Springer Science \& Business Media},
  \bibinfo{year}{2012}.
\bibitem[{Brezzi et~al.(2004)Brezzi, Marini, and S{\"u}li}]{brezzi2004}
\bibinfo{author}{F.~Brezzi}, \bibinfo{author}{L.~D. Marini},
  \bibinfo{author}{E.~S{\"u}li},
\newblock \bibinfo{title}{{Discontinuous Galerkin methods for first-order
  hyperbolic problems}},
\newblock \bibinfo{journal}{{Mathematical Models and Methods in Applied
  Sciences}} \bibinfo{volume}{14} (\bibinfo{year}{2004})
  \bibinfo{pages}{1893--1903}.
\bibitem[{Ern and Guermond(2006)}]{ern2006}
\bibinfo{author}{A.~Ern}, \bibinfo{author}{J.-L. Guermond},
\newblock \bibinfo{title}{{Discontinuous Galerkin methods for Friedrichs'
  systems. I. General theory}},
\newblock \bibinfo{journal}{{SIAM Journal on Numerical Analysis}}
  \bibinfo{volume}{44} (\bibinfo{year}{2006}) \bibinfo{pages}{753--778}.
\bibitem[{Bochev and Gunzburger(2009)}]{bochev2009}
\bibinfo{author}{P.~B. Bochev}, \bibinfo{author}{M.~D. Gunzburger},
  \bibinfo{title}{Least-squares finite element methods}, volume
  \bibinfo{volume}{166}, \bibinfo{publisher}{Springer Science \& Business
  Media}, \bibinfo{year}{2009}.
\bibitem[{Demkowicz and Gopalakrishnan(2010)}]{demkowicz2010}
\bibinfo{author}{L.~Demkowicz}, \bibinfo{author}{J.~Gopalakrishnan},
\newblock \bibinfo{title}{{A class of discontinuous Petrov--Galerkin methods.
  Part I: The transport equation}},
\newblock \bibinfo{journal}{{Computer Methods in Applied Mechanics and
  Engineering}} \bibinfo{volume}{199} (\bibinfo{year}{2010})
  \bibinfo{pages}{1558--1572}.
\bibitem[{Demkowicz and Gopalakrishnan(2011)}]{demkowicz2011}
\bibinfo{author}{L.~Demkowicz}, \bibinfo{author}{J.~Gopalakrishnan},
\newblock \bibinfo{title}{{A class of discontinuous Petrov--Galerkin methods.
  II. Optimal test functions}},
\newblock \bibinfo{journal}{{Numerical Methods for Partial Differential
  Equations}} \bibinfo{volume}{27} (\bibinfo{year}{2011})
  \bibinfo{pages}{70--105}.
\bibitem[{Demkowicz and Gopalakrishnan(2014)}]{DemGopBOOK-CH2014}
\bibinfo{author}{L.~Demkowicz}, \bibinfo{author}{J.~Gopalakrishnan},
\newblock \bibinfo{title}{An overview of the discontinuous {Petrov Galerkin}
  method},
\newblock in: \bibinfo{editor}{X.~Feng}, \bibinfo{editor}{O.~Karakashian},
  \bibinfo{editor}{Y.~Xing} (Eds.), \bibinfo{booktitle}{Recent Developments in
  Discontinuous Galerkin Finite Element Methods for Partial Differential
  Equations: 2012 John H Barrett Memorial Lectures}, volume
  \bibinfo{volume}{157} of \textit{\bibinfo{series}{The IMA Volumes in
  Mathematics and its Applications}}, \bibinfo{publisher}{Springer},
  \bibinfo{address}{Cham}, \bibinfo{year}{2014}, pp. \bibinfo{pages}{149--180}.
\bibitem[{Calo et~al.(2014)Calo, Collier, and Niemi}]{Calo2014dPG}
\bibinfo{author}{V.~M. Calo}, \bibinfo{author}{N.~O. Collier},
  \bibinfo{author}{A.~H. Niemi},
\newblock \bibinfo{title}{{Analysis of the discontinuous Petrov–Galerkin
  method with optimal test functions for the Reissner–Mindlin plate bending
  model}},
\newblock \bibinfo{journal}{Computers \& Mathematics with Applications}
  \bibinfo{volume}{66} (\bibinfo{year}{2014}) \bibinfo{pages}{2570 -- 2586}.
\bibitem[{Demkowicz et~al.(2012)Demkowicz, Gopalakrishnan, and
  Niemi}]{Demkowicz2012}
\bibinfo{author}{L.~Demkowicz}, \bibinfo{author}{J.~Gopalakrishnan},
  \bibinfo{author}{A.~H. Niemi},
\newblock \bibinfo{title}{{A class of discontinuous Petrov–Galerkin methods.
  Part III: Adaptivity}},
\newblock \bibinfo{journal}{Applied Numerical Mathematics} \bibinfo{volume}{62}
  (\bibinfo{year}{2012}) \bibinfo{pages}{396 -- 427}. \bibinfo{note}{Third
  Chilean Workshop on Numerical Analysis of Partial Differential Equations
  (WONAPDE 2010)}.
\bibitem[{Niemi et~al.(2011)Niemi, Bramwell, and Demkowicz}]{Niemi2011a}
\bibinfo{author}{A.~Niemi}, \bibinfo{author}{J.~Bramwell},
  \bibinfo{author}{L.~Demkowicz},
\newblock \bibinfo{title}{{Discontinuous Petrov-Galerkin method with optimal
  test functions for thin-body problems in solid mechanics}},
\newblock \bibinfo{journal}{Computer Methods in Applied Mechanics and
  Engineering} \bibinfo{volume}{200} (\bibinfo{year}{2011})
  \bibinfo{pages}{1291--1300}.
\bibitem[{Niemi et~al.(2013)Niemi, Collier, and Calo}]{Niemi2013}
\bibinfo{author}{A.~H. Niemi}, \bibinfo{author}{N.~O. Collier},
  \bibinfo{author}{V.~M. Calo},
\newblock \bibinfo{title}{{Automatically stable discontinuous Petrov–Galerkin
  methods for stationary transport problems: Quasi-optimal test space norm}},
\newblock \bibinfo{journal}{Computers \& Mathematics with Applications}
  \bibinfo{volume}{66} (\bibinfo{year}{2013}) \bibinfo{pages}{2096 -- 2113}.
  \bibinfo{note}{ICNC-FSKD 2012}.
\bibitem[{Niemi et~al.(2011)Niemi, Collier, and Calo}]{Niemi2011}
\bibinfo{author}{A.~H. Niemi}, \bibinfo{author}{N.~O. Collier},
  \bibinfo{author}{V.~M. Calo},
\newblock \bibinfo{title}{{Discontinuous Petrov-Galerkin method based on the
  optimal test space norm for one-dimensional transport problems}},
\newblock \bibinfo{journal}{Procedia Computer Science} \bibinfo{volume}{4}
  (\bibinfo{year}{2011}) \bibinfo{pages}{1862 -- 1869}.
  \bibinfo{note}{Proceedings of the International Conference on Computational
  Science, ICCS 2011}.
\bibitem[{Cier et~al.(2020{\natexlab{a}})Cier, Poulet, Rojas, Calo, and
  Veveakis}]{cier2020}
\bibinfo{author}{R.~J. Cier}, \bibinfo{author}{T.~Poulet},
  \bibinfo{author}{S.~Rojas}, \bibinfo{author}{V.~M. Calo},
  \bibinfo{author}{M.~Veveakis},
\newblock \bibinfo{title}{Adaptive stabilized finite elements: Continuation
  analysis of compaction banding in geomaterials},
\newblock \bibinfo{journal}{arXiv preprint arXiv:2008.01396}
  (\bibinfo{year}{2020}{\natexlab{a}}).
\bibitem[{Cier et~al.(2020{\natexlab{b}})Cier, Rojas, and Calo}]{cier2020_2}
\bibinfo{author}{R.~J. Cier}, \bibinfo{author}{S.~Rojas},
  \bibinfo{author}{V.~M. Calo},
\newblock \bibinfo{title}{A nonlinear weak constraint enforcement method for
  advection-dominated diffusion problems},
\newblock \bibinfo{journal}{Mechanics Research Communications}
  (\bibinfo{year}{2020}{\natexlab{b}}) \bibinfo{pages}{103602}.
\bibitem[{{\L}o{\'s} et~al.(2020){\L}o{\'s}, Rojas, Paszy{\'n}ski, Muga, and
  Calo}]{los2020}
\bibinfo{author}{M.~{\L}o{\'s}}, \bibinfo{author}{S.~Rojas},
  \bibinfo{author}{M.~Paszy{\'n}ski}, \bibinfo{author}{I.~Muga},
  \bibinfo{author}{V.~M. Calo},
\newblock \bibinfo{title}{{A stable discontinuous Galerkin based isogeometric
  residual minimization for the Stokes problem}},
\newblock in: \bibinfo{booktitle}{International Conference on Computational
  Science}, \bibinfo{organization}{Springer}, \bibinfo{year}{2020}, pp.
  \bibinfo{pages}{197--211}.
\bibitem[{Kyburg et~al.(2020)Kyburg, Rojas, and Calo}]{kyburg2020}
\bibinfo{author}{F.~Kyburg}, \bibinfo{author}{S.~Rojas}, \bibinfo{author}{V.~M.
  Calo},
\newblock \bibinfo{title}{Incompressible flow modeling using an adaptive
  stabilized finite element method based on residual minimization},
\newblock \bibinfo{journal}{arXiv preprint arXiv:2011.09182}
  (\bibinfo{year}{2020}).
\bibitem[{Rojas et~al.(2020)Rojas, Pardo, Behnoudfar, and Calo}]{rojas2020}
\bibinfo{author}{S.~Rojas}, \bibinfo{author}{D.~Pardo},
  \bibinfo{author}{P.~Behnoudfar}, \bibinfo{author}{V.~M. Calo},
\newblock \bibinfo{title}{Residual minimization for goal-oriented adaptivity},
\newblock \bibinfo{journal}{arXiv preprint arXiv:2007.08824}
  (\bibinfo{year}{2020}).
\bibitem[{Codina(1998)}]{codina1998}
\bibinfo{author}{R.~Codina},
\newblock \bibinfo{title}{Comparison of some finite element methods for solving
  the diffusion-convection-reaction equation},
\newblock \bibinfo{journal}{Computer methods in applied mechanics and
  engineering} \bibinfo{volume}{156} (\bibinfo{year}{1998})
  \bibinfo{pages}{185--210}.
\bibitem[{Hughes et~al.(2018)Hughes, Scovazzi, and Franca}]{hughes2018}
\bibinfo{author}{T.~J. Hughes}, \bibinfo{author}{G.~Scovazzi},
  \bibinfo{author}{L.~P. Franca},
\newblock \bibinfo{title}{Multiscale and stabilized methods},
\newblock \bibinfo{journal}{Encyclopedia of Computational Mechanics Second
  Edition}  (\bibinfo{year}{2018}) \bibinfo{pages}{1--64}.
\bibitem[{Di~Pietro et~al.(2008)Di~Pietro, Ern, and Guermond}]{di2008}
\bibinfo{author}{D.~A. Di~Pietro}, \bibinfo{author}{A.~Ern},
  \bibinfo{author}{J.-L. Guermond},
\newblock \bibinfo{title}{{Discontinuous Galerkin methods for anisotropic
  semidefinite diffusion with advection}},
\newblock \bibinfo{journal}{SIAM Journal on Numerical Analysis}
  \bibinfo{volume}{46} (\bibinfo{year}{2008}) \bibinfo{pages}{805--831}.
\bibitem[{Di~Pietro and Ern(2012)}]{di2011mathematical}
\bibinfo{author}{D.~A. Di~Pietro}, \bibinfo{author}{A.~Ern},
  \bibinfo{title}{{Mathematical aspects of discontinuous Galerkin methods}},
  volume~\bibinfo{volume}{69}, \bibinfo{publisher}{Springer Science},
  \bibinfo{year}{2012}.
\bibitem[{Shahbazi(2005)}]{shahbazi2005}
\bibinfo{author}{K.~Shahbazi},
\newblock \bibinfo{title}{An explicit expression for the penalty parameter of
  the interior penalty method},
\newblock \bibinfo{journal}{Journal of Computational Physics}
  \bibinfo{volume}{205} (\bibinfo{year}{2005}) \bibinfo{pages}{401--407}.
\bibitem[{Epshteyn and Rivi{\`e}re(2007)}]{epshteyn2007}
\bibinfo{author}{Y.~Epshteyn}, \bibinfo{author}{B.~Rivi{\`e}re},
\newblock \bibinfo{title}{{Estimation of penalty parameters for symmetric
  interior penalty Galerkin methods}},
\newblock \bibinfo{journal}{Journal of Computational and Applied Mathematics}
  \bibinfo{volume}{206} (\bibinfo{year}{2007}) \bibinfo{pages}{843--872}.
\bibitem[{Hartmann and Houston(2008)}]{hartmann2008}
\bibinfo{author}{R.~Hartmann}, \bibinfo{author}{P.~Houston},
\newblock \bibinfo{title}{{An optimal order interior penalty discontinuous
  Galerkin discretization of the compressible Navier--Stokes equations}},
\newblock \bibinfo{journal}{Journal of Computational Physics}
  \bibinfo{volume}{227} (\bibinfo{year}{2008}) \bibinfo{pages}{9670--9685}.
\bibitem[{Bastian et~al.(2012)Bastian, Blatt, and Scheichl}]{bastian2012}
\bibinfo{author}{P.~Bastian}, \bibinfo{author}{M.~Blatt},
  \bibinfo{author}{R.~Scheichl},
\newblock \bibinfo{title}{Algebraic multigrid for discontinuous galerkin
  discretizations of heterogeneous elliptic problems},
\newblock \bibinfo{journal}{Numerical Linear Algebra with Applications}
  \bibinfo{volume}{19} (\bibinfo{year}{2012}) \bibinfo{pages}{367--388}.
\bibitem[{Wheeler(1978)}]{wheeler1978}
\bibinfo{author}{M.~F. Wheeler},
\newblock \bibinfo{title}{An elliptic collocation-finite element method with
  interior penalties},
\newblock \bibinfo{journal}{SIAM Journal on Numerical Analysis}
  \bibinfo{volume}{15} (\bibinfo{year}{1978}) \bibinfo{pages}{152--161}.
\bibitem[{Arnold(1982)}]{arnold1982}
\bibinfo{author}{D.~N. Arnold},
\newblock \bibinfo{title}{An interior penalty finite element method with
  discontinuous elements},
\newblock \bibinfo{journal}{SIAM Journal on Numerical Analysis}
  \bibinfo{volume}{19} (\bibinfo{year}{1982}) \bibinfo{pages}{742--760}.
\bibitem[{Arnold et~al.(2002)Arnold, Brezzi, Cockburn, and
  Marini}]{arnold2002unified}
\bibinfo{author}{D.~N. Arnold}, \bibinfo{author}{F.~Brezzi},
  \bibinfo{author}{B.~Cockburn}, \bibinfo{author}{L.~D. Marini},
\newblock \bibinfo{title}{{Unified analysis of discontinuous Galerkin methods
  for elliptic problems}},
\newblock \bibinfo{journal}{SIAM Journal on Numerical Analysis}
  \bibinfo{volume}{39} (\bibinfo{year}{2002}) \bibinfo{pages}{1749--1779}.
\bibitem[{Cohen et~al.(2012)Cohen, Dahmen, and Welper}]{cohen2012}
\bibinfo{author}{A.~Cohen}, \bibinfo{author}{W.~Dahmen},
  \bibinfo{author}{G.~Welper},
\newblock \bibinfo{title}{Adaptivity and variational stabilization for
  convection-diffusion equations},
\newblock \bibinfo{journal}{ESAIM: Mathematical Modelling and Numerical
  Analysis} \bibinfo{volume}{46} (\bibinfo{year}{2012})
  \bibinfo{pages}{1247--1273}.
\bibitem[{Karakashian and Pascal(2003)}]{karakashian2003}
\bibinfo{author}{O.~A. Karakashian}, \bibinfo{author}{F.~Pascal},
\newblock \bibinfo{title}{A posteriori error estimates for a discontinuous
  galerkin approximation of second-order elliptic problems},
\newblock \bibinfo{journal}{SIAM Journal on Numerical Analysis}
  \bibinfo{volume}{41} (\bibinfo{year}{2003}) \bibinfo{pages}{2374--2399}.
\bibitem[{Burman and Ern(2007)}]{burman2007}
\bibinfo{author}{E.~Burman}, \bibinfo{author}{A.~Ern},
\newblock \bibinfo{title}{Continuous interior penalty hp-finite element methods
  for advection and advection-diffusion equations},
\newblock \bibinfo{journal}{Mathematics of computation} \bibinfo{volume}{76}
  (\bibinfo{year}{2007}) \bibinfo{pages}{1119--1140}.
\bibitem[{Ern and Guermond(2017)}]{ern2017}
\bibinfo{author}{A.~Ern}, \bibinfo{author}{J.-L. Guermond},
\newblock \bibinfo{title}{Finite element quasi-interpolation and best
  approximation},
\newblock \bibinfo{journal}{ESAIM: Mathematical Modelling and Numerical
  Analysis} \bibinfo{volume}{51} (\bibinfo{year}{2017})
  \bibinfo{pages}{1367--1385}.
\bibitem[{Aln{\ae}s et~al.(2015)Aln{\ae}s, Blechta, Hake, Johansson, Kehlet,
  Logg, Richardson, Ring, Rognes, and Wells}]{alnaes2015fenics}
\bibinfo{author}{M.~S. Aln{\ae}s}, \bibinfo{author}{J.~Blechta},
  \bibinfo{author}{J.~Hake}, \bibinfo{author}{A.~Johansson},
  \bibinfo{author}{B.~Kehlet}, \bibinfo{author}{A.~Logg},
  \bibinfo{author}{C.~Richardson}, \bibinfo{author}{J.~Ring},
  \bibinfo{author}{M.~E. Rognes}, \bibinfo{author}{G.~N. Wells},
\newblock \bibinfo{title}{{The FEniCS project version 1.5}},
\newblock \bibinfo{journal}{Archive of Numerical Software} \bibinfo{volume}{3}
  (\bibinfo{year}{2015}) \bibinfo{pages}{9--23}.
\bibitem[{D{\"o}rfler(1996)}]{dorfler1996convergent}
\bibinfo{author}{W.~D{\"o}rfler},
\newblock \bibinfo{title}{A convergent adaptive algorithm for {P}oisson's
  equation},
\newblock \bibinfo{journal}{SIAM Journal on Numerical Analysis}
  \bibinfo{volume}{33} (\bibinfo{year}{1996}) \bibinfo{pages}{1106--1124}.
\bibitem[{Bank et~al.(1983)Bank, Sherman, and Weiser}]{bank1983some}
\bibinfo{author}{R.~E. Bank}, \bibinfo{author}{A.~H. Sherman},
  \bibinfo{author}{A.~Weiser},
\newblock \bibinfo{title}{Some refinement algorithms and data structures for
  regular local mesh refinement},
\newblock \bibinfo{journal}{Scientific Computing, Applications of Mathematics
  and Computing to the Physical Sciences} \bibinfo{volume}{1}
  (\bibinfo{year}{1983}) \bibinfo{pages}{3--17}.
\bibitem[{Bank et~al.(1989)Bank, Welfert, and Yserentant}]{bank1989class}
\bibinfo{author}{R.~E. Bank}, \bibinfo{author}{B.~D. Welfert},
  \bibinfo{author}{H.~Yserentant},
\newblock \bibinfo{title}{A class of iterative methods for solving saddle point
  problems},
\newblock \bibinfo{journal}{Numerische Mathematik} \bibinfo{volume}{56}
  (\bibinfo{year}{1989}) \bibinfo{pages}{645--666}.
\bibitem[{Chen et~al.(2008)Chen, Davis, Hager, and
  Rajamanickam}]{chen2008algorithm}
\bibinfo{author}{Y.~Chen}, \bibinfo{author}{T.~A. Davis},
  \bibinfo{author}{W.~W. Hager}, \bibinfo{author}{S.~Rajamanickam},
\newblock \bibinfo{title}{{Algorithm 887: CHOLMOD, supernodal sparse Cholesky
  factorization and update/downdate}},
\newblock \bibinfo{journal}{ACM Transactions on Mathematical Software (TOMS)}
  \bibinfo{volume}{35} (\bibinfo{year}{2008}) \bibinfo{pages}{22}.
\bibitem[{Mitchell(2013)}]{mitchell2013collection}
\bibinfo{author}{W.~F. Mitchell},
\newblock \bibinfo{title}{{A collection of 2D elliptic problems for testing
  adaptive grid refinement algorithms}},
\newblock \bibinfo{journal}{Applied Mathematics and Computation}
  \bibinfo{volume}{220} (\bibinfo{year}{2013}) \bibinfo{pages}{350--364}.
\bibitem[{Oden and Patra(1995)}]{oden1995parallel}
\bibinfo{author}{J.~T. Oden}, \bibinfo{author}{A.~Patra},
\newblock \bibinfo{title}{A parallel adaptive strategy for hp finite element
  computations},
\newblock \bibinfo{journal}{Computer Methods in Applied Mechanics and
  Engineering} \bibinfo{volume}{121} (\bibinfo{year}{1995})
  \bibinfo{pages}{449--470}.
\bibitem[{Burman and Zunino(2006)}]{burman2006domain}
\bibinfo{author}{E.~Burman}, \bibinfo{author}{P.~Zunino},
\newblock \bibinfo{title}{A domain decomposition method based on weighted
  interior penalties for advection-diffusion-reaction problems},
\newblock \bibinfo{journal}{SIAM Journal on Numerical Analysis}
  \bibinfo{volume}{44} (\bibinfo{year}{2006}) \bibinfo{pages}{1612--1638}.

\end{thebibliography}

\end{document}